\documentclass[10pt]{article}

\usepackage{geometry}
\usepackage{setspace}
\usepackage{graphicx}
\usepackage{multirow}
\usepackage{amsmath,amssymb,amsfonts}
\usepackage{amsthm}
\usepackage{mathrsfs}
\usepackage{xcolor}
\usepackage{booktabs}
\usepackage{algorithm}
\usepackage{algorithmicx}
\usepackage{algpseudocode}
\usepackage{listings}
\usepackage{longtable}
\usepackage{mathtools}
\usepackage{natbib}
\usepackage{arydshln}
\usepackage{hyperref}
\usepackage{animate}
\usepackage{caption}
\usepackage{subcaption}
\usepackage{bm}
\usepackage{orcidlink}
\usepackage{cleveref}
\usepackage{xr}

\newtheorem{condition}{Condition}

\usepackage{authblk}

\hypersetup{colorlinks=true, linkcolor=blue, citecolor=blue}

\usepackage{bm}

\newenvironment{biography}[1][]{\par\medskip\noindent\textbf{#1}\par\medskip\noindent}{\par\medskip}

\title{A Partially Functional Dynamic Structural Equation Model for Multi-Resolution Environmental Data}

\author[1]{Xinwen Liu}
\author[1]{Niansheng Tang}
\author[2]{Song Xi Chen}

\affil[1]{Yunnan Key Laboratory of Statistical Modeling and Data Analysis; Yunnan International Joint Laboratory for Carbon Emission Statistics Accounting and AI-driven Forecasting \& Warning, Yunnan University, ChengGon, KunMing 650500, Yunnan, China}
\affil[2]{Department of Statistics and Data Science, Tsinghua University, Haidian, Beijing 100084, China}

\begin{document}

\maketitle

\begin{abstract}
Understanding the complex relationships between atmospheric pollutant emissions and their multifaceted determinants presents a dual challenge: driving factors operate at fundamentally different temporal resolutions, from continuously monitored meteorological variables to annually reported socio-economic indicators, and their interconnections evolve dynamically over time. To address these challenges, we propose a Partially Functional Dynamic Structural Equation Model (PFDSEM) that coherently integrates functional covariates (e.g., high-frequency meteorological data) and scalar predictors (e.g., economic and demographic indicators) within a unified dynamic structural framework. The model captures non-stationary temporal dependencies and inter-variable correlations via a Conditional Autoregressive (CAR) structure combined with a Linear Model of Coregionalization (LMC), while functional covariates are incorporated through basis expansion with Bayesian P-spline smoothing. Bayesian inference via Markov Chain Monte Carlo provides full uncertainty quantification. Comprehensive simulation studies confirm accurate parameter recovery and robustness to prior specifications under diverse conditions. Applying the PFDSEM to pollutant emissions data from 30 Chinese provinces (2015--2020), we identify temporally dynamic and province-specific associations between ten categories of socio-environmental factors and ten major air pollutants, including CO$_2$. The results reveal substantial cross-province heterogeneity in the strength and direction of these associations and pronounced seasonal patterns in meteorological effects, offering a quantitative evidence base for designing temporally adaptive and regionally tailored environmental policies.
\end{abstract}


\section{Introduction}\label{sec:Intro.}

\begin{figure}[!ht]
\centering
\includegraphics[width=1\textwidth]{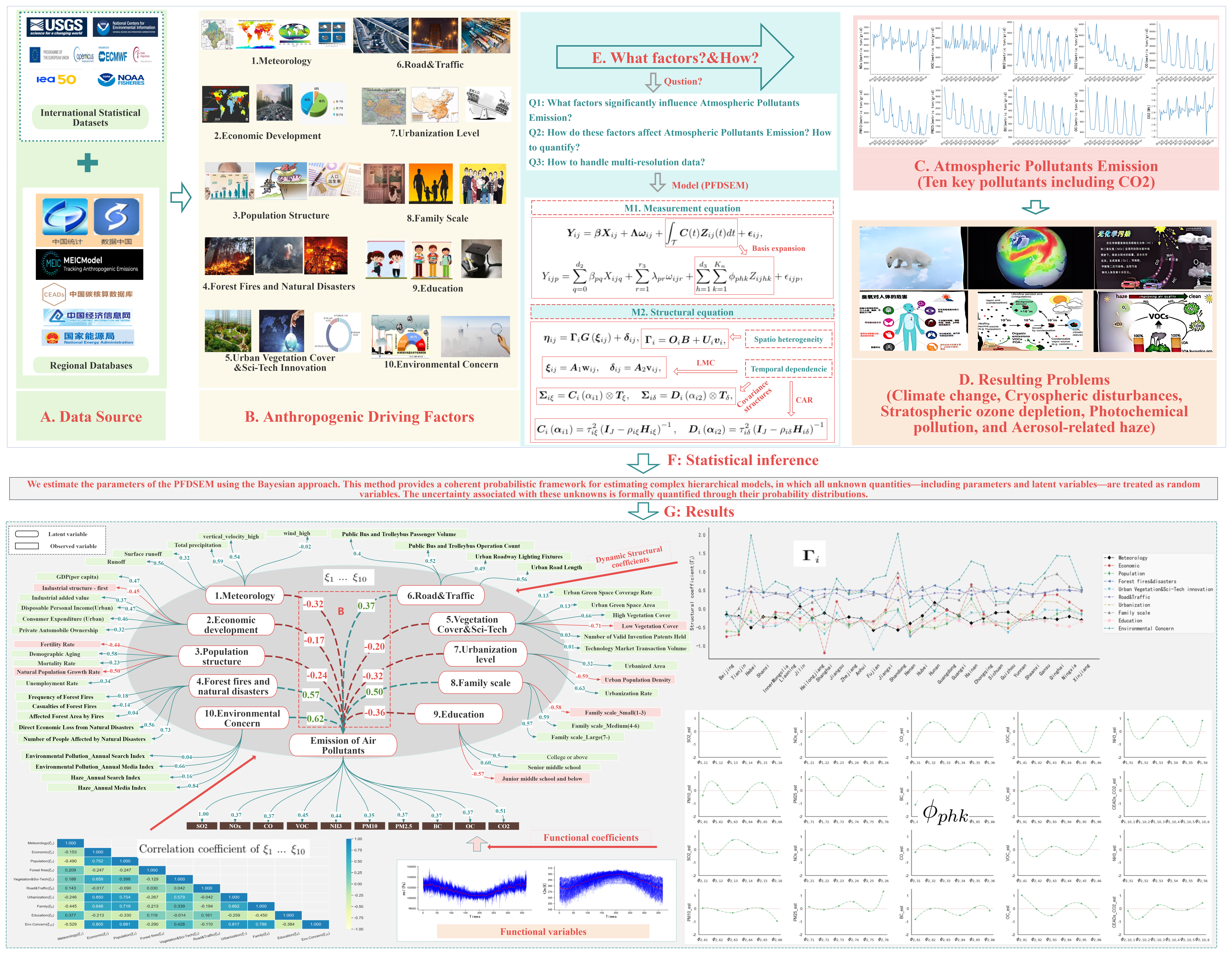}
\caption{\textbf{Panel A}: Illustration of multi-resolution data sources: functional meteorological data (e.g., hourly temperature, pressure) and scalar socio-economic indicators (e.g., annual GDP, population). 
\textbf{Panel B}: A taxonomy of ten socio-environmental factor categories potentially associated with atmospheric pollutant emissions, grouped to illustrate interconnected domains. 
\textbf{Panel C}: Schematic diagram illustrating temporal trends and seasonal patterns in emissions of ten key pollutants (including CO$_2$) across the study period. 
\textbf{Panel D}: The primary environmental consequences of pollutant emissions: Climate Change, Cryospheric Disturbances, Stratospheric Ozone Depletion, Photochemical Pollution, and Aerosol-related Haze. 
\textbf{Panel E}: The PFDSEM framework addresses three core analytical questions: (Q1) Which factors are most strongly associated with emissions? (Q2) How do these associations vary across provinces and over time? (Q3) How to coherently integrate data measured at different temporal resolutions? 
\textbf{Panel F}: Bayesian inference for the PFDSEM, where all unknown quantities are treated as random variables and uncertainty is quantified through posterior probability distributions. 
\textbf{Panel G}: Key outputs include posterior estimates of factor loadings, structural coefficients, temporal dependence parameters (CAR and LMC), and functional coefficients, enabling comprehensive quantification of the dynamic, multi-resolution associations between environmental factors and pollutant emissions.}
\label{fig:PFDSEM}
\end{figure}

Numerous climate and environmental issue trajectories, including those of anthropogenic climate change, cryosphere perturbations, stratospheric ozone depletion, photochemical pollution, and particulate matter-related smog, are intimately connected with human activities that release various greenhouse gases and pollutants into the Earth's atmosphere 
 (Fig. \ref{fig:PFDSEM} Panel D). In the field of atmospheric science and climate change, a significant challenge is presented by the need to accurately define the scale, spatial and temporal distribution, representation, and prediction of atmospheric pollutant emissions, as well as to explore the driving role of human factors in these emission changes (Fig. \ref{fig:PFDSEM} Panel B and C). These challenges are in lockstep with the strategic approach highlighted at the 20th National Congress of the Communist Party of China, which underscores the importance of a dual strategy focused on ``promoting ecological priorities, conservation, efficiency, and green, low-carbon development''  \citep{xi2023china}. Therefore, comprehensive research and evaluation of the spatial and temporal patterns of carbon emission indicators, the underlying driving factors, and their interconnections are essential to develop scientifically robust and accurate emission reduction strategies (Fig. \ref{fig:PFDSEM} Panel E).

A review of the extant literature reveals that traditional regression frameworks and spatial econometrics methods are limited in their capacity to incorporate dynamic temporal interdependencies between variables \citep{apergis2010renewable,ghosh2010examining,jayanthakumaran2012co2,ozcan2013nexus}. These models predominantly investigate unidirectional relationships, typically from multiple independent variables to a single dependent variable, whereas the complexity of environmental systems necessitates simultaneous examination of multiple interrelated response and explanatory variables \citep{chaudhary2018factors,rauf2018structural,song2018environmental}. 
The complexity of variable relationships in empirical studies often exceeds this simplistic representation, necessitating a simultaneous examination of multiple interrelated response and explanatory variables.

Structural equation modeling (SEM) offers a well-established multivariate framework for analyzing complex relationships between observed and latent variables. However, conventional SEM formulations typically ignore temporal dynamics in observed indicators. Addressing this gap, \citet{song2012bayesian} proposed the Generalized Random Coefficient Structural Equation Model (GRCSEM), which explicitly incorporates temporal correlations and cross-sectional heterogeneity. This framework significantly advances the analysis of complex, dynamic variable interactions. Yet, the GRCSEM is designed for scalar variables and cannot natively accommodate functional covariates measured at higher temporal resolutions.

The analysis of functional data, characterized by observations indexed over a continuum (e.g., time, wavelength, spatial location), has attracted substantial methodological attention, particularly in environmental sciences \citep{Ramsay2005FDA}. In atmospheric research, functional time-series, such as hourly temperature, humidity, and wind vector profiles, support meteorological forecasting, while spatially indexed pollutant measurements enable air-quality monitoring. Functional data analysis has a rich tradition in environmental applications, from modeling air pollution diurnal cycles to predicting precipitation patterns. However, integrating functional predictors into a dynamic structural equation framework with latent variables and temporal dependence remains an open methodological challenge.

To address these gaps, we propose a Partially Functional Dynamic Structural Equation Model (PFDSEM) for analyzing complex temporal and cross-sectional structural relationships among multivariate indicators measured at mixed resolutions. Unlike standard spatio-temporal panel models, which are typically limited to univariate or multivariate scalar responses, our framework offers three distinct advantages: 
\begin{enumerate}
    \item[i)] \textit{Multi-resolution data integration.} We employ basis expansion techniques to project infinite-dimensional functional covariates onto a finite-dimensional space spanned by a set of basis functions (e.g., Fourier, B-spline, or wavelet bases). This transforms the functional regression problem into a computationally tractable form without sacrificing the continuous nature of the data. Unlike ad-hoc discretization or two-step smoothing approaches, our method preserves functional continuity and enables integrated estimation within the full model, ensuring that uncertainty from the functional representation propagates coherently through all levels of inference.    
    \item[ii)] \textit{Dynamic temporal structure.} Our model incorporates temporal dynamics through a Conditional Autoregressive (CAR) structure combined with a Linear Model of Coregionalization (LMC), explicitly accounting for both inter-variable and inter-time dependency structures. This allows for effective modeling of non-stationary temporal correlations common in longitudinal environmental studies. The framework also accommodates province-level heterogeneity in structural parameters via random effects and supports nonlinear relationships among latent variables.
        \item[iii)] \textit{Flexible functional modeling with uncertainty quantification.} We adopt Bayesian P-splines for the nonparametric approximation of coefficient functions, ensuring smoothness and guarding against overfitting through a hierarchical prior structure on the smoothing parameters. The fully Bayesian approach provides posterior uncertainty bands for the estimated functional coefficients, enabling formal statistical assessment of seasonal and sub-seasonal patterns.
\end{enumerate}

This synergistic combination of a flexible functional representation, a dynamic latent variable structure, and a fully Bayesian inferential framework allows the PFDSEM to handle complex, multi-resolution problems that are computationally prohibitive for conventional estimation techniques. Simulation studies demonstrate accurate parameter recovery across varying sample sizes, functional complexities, and prior specifications. In an empirical application analyzing China's provincial atmospheric pollutant emission indicators (APEI) and their determinants (2015--2020), the PFDSEM reveals substantial province-level heterogeneity and pronounced seasonal patterns in the strength of key associations, offering a nuanced evidence base for designing temporally adaptive and regionally differentiated environmental policies  (Fig. \ref{fig:PFDSEM} Panel G).

The environmental-scientific contributions of this paper are three-fold. First, we provide the first unified quantification of how ten categories of anthropogenic and natural factors are jointly associated with ten major air pollutants across 30 Chinese provinces, respecting the intrinsic multi-resolution nature of the underlying data. Second, we reveal substantial province-level heterogeneity in these associations, challenging the one-size-fits-all assumption common in national-level policy analysis and highlighting the need for regionally tailored interventions. Third, by explicitly modeling temporal non-stationarity, we identify seasonal and inter-annual shifts in the strength of key associations, offering evidence for temporally adaptive environmental regulation.

Our paper is organized as follows. Section \ref{sec:PFDSEM} presents the full specification of the PFDSEM, including the measurement equation with functional and scalar covariates, the structural equation with province-level heterogeneity, and the temporal dependence structure captured via the LMC combined with CAR priors (Section \ref{sec:MS}). We provide a detailed discussion of model identifiability, extending standard conditions to cover the LMC parameters (Section \ref{sec:MIC}). Section \ref{sec:BE} develops the Bayesian inference framework, outlining the prior specifications (Section \ref{sec:PS}), deriving the full conditional distributions for all parameters, and describing the hybrid MCMC algorithm used for posterior computation (Section \ref{sec:CD} and Appendices). Section \ref{sec:SS} reports comprehensive simulation studies assessing parameter recovery and prior sensitivity.  Section \ref{sec:Application} applies the PFDSEM to analyze atmospheric pollutant emissions across 30 Chinese provinces from 2015--2020. Section \ref{sec:Disc} provides a discussion of the findings, policy implications, methodological limitations, and directions for future research. Technical derivations, additional simulation results, and supplementary empirical tables and figures are collected in the Appendices.

\section{Partially Functional Dynamic Structural Equation Model}\label{sec:PFDSEM} 

\subsection{Data Characteristics and Modeling Challenges}\label{sec:DCMC}

The analysis of atmospheric pollutant emissions and their determinants presents several interconnected data challenges that motivate the development of the PFDSEM. We articulate these challenges explicitly to clarify how each model component addresses a specific aspect of the data structure.

\begin{description}
    \item[Challenge 1: Multi-resolution measurements.] The covariates hypothesized to be associated with emissions operate at fundamentally different temporal granularities. Meteorological variables (e.g., temperature, atmospheric pressure) are continuously monitored and are best represented as smooth functions of time within each observation period. Socio-economic indicators (e.g., GDP, urbanization rate) are discrete scalar values reported at annual or quarterly frequencies. A principled modeling framework must respect these ontological differences without resorting to ad-hoc aggregation of high-frequency data (which discards information about within-period variability) or artificial discretization of functional data (which ignores smoothness and continuity). \textit{PFDSEM solution: Basis expansion for functional covariates (Section \ref{sec:MS}, Eq. \eqref{eq: PFSEM_reduced}).} 
    
    \item[Challenge 2: Latent confounding and multiple interrelated outcomes.] Ten distinct pollutant species are measured simultaneously, and they share common anthropogenic and natural sources. Individual regression models for each pollutant ignore the rich covariance structure among pollutants and among driving factors, potentially leading to inefficient inference and inability to model synergistic effects. Moreover, many driving factors (e.g., ``economic development,'' ``technological innovation'') are themselves latent constructs measured through multiple observed indicators, requiring a measurement model. \textit{PFDSEM solution: Latent variable measurement model with structured factor loadings (Section \ref{sec:MS}, Eq. \eqref{eq: PFSEM}).}
    
    \item[Challenge 3: Temporal dynamics and non-stationarity.] The relationships between emissions and their predictors are unlikely to be static over time. Policy interventions, economic cycles, and gradual technological change can alter the strength and even the sign of associations across years or seasons. Moreover, observations at adjacent time points are correlated, violating the independence assumptions of standard regression models and requiring explicit modeling of temporal dependence. \textit{PFDSEM solution: CAR(1) temporal structure with LMC for inter-variable covariances (Section \ref{sec:MS}, Eqs. \eqref{eq: covariance_matrices}--\eqref{eq: CAR}).}
    
    \item[Challenge 4: Province-level heterogeneity.] China's provinces vary dramatically in industrial structure, energy mix, demographic profile, and environmental policy implementation. A model assuming homogeneous structural parameters across all provinces would mask this critical heterogeneity and produce misleading aggregate conclusions. \textit{PFDSEM solution: Random coefficients for structural parameters (Section \ref{sec:MS}, Eq. \eqref{eq: structural_heterogeneity}).}
\end{description}

By explicitly recognizing these four challenges, the model specification that follows can be understood as a coherent, purpose-built solution rather than an arbitrary assembly of statistical components. Each element of the PFDSEM is motivated by a specific, well-characterized feature of the environmental data.

\subsection{Model and Notation}\label{sec:MS} 
 We consider $n$ individuals observed at $J$ time points. For individual $i=1, \ldots, n$ and time $j=1,\ldots, J$, let $\boldsymbol{Y}_{i j}$ be a $d_1 \times 1$ vector of manifest variables, $\boldsymbol{X}_{i j}=(1, X_{i j 1},\cdots,X_{i j d_2})$  represent a $(d_2+1) \times 1$ vector of scalar covariates, $\boldsymbol{Z}_{i j}(t)=(Z_{i j 1}(t),\cdots,Z_{i j d_3}(t))^\top$ be $d_3 \times 1$ functional covariates, where $t \in [0,T]$ with $T$ denoting the period length. The latent vector $\boldsymbol{\omega}_{ij} = (\boldsymbol{\eta}_{ij}^\top, \boldsymbol{\xi}_{ij}^\top)^\top$ partitions into outcome latent variables $\boldsymbol{\eta}_{ij} \, (r_1 \times 1)$ and explanatory latent variables $\boldsymbol{\xi}_{ij} \, (r_2 \times 1)$, with $r_1 + r_2 = r_3$. 

The manifest variables in $\boldsymbol{Y}_{i j}$ are assumed to relate to observed scalar covariates $\boldsymbol{X}_{i j}$, functional covariates $\boldsymbol{Z}_{i j}(t)$,  and unobserved latent variables $\boldsymbol{\omega}_{i j}$ through a partially functional dynamic structural equation model 
\begin{equation}\label{eq: PFSEM}
\begin{aligned} 
&\boldsymbol{Y}_{i j}=\boldsymbol{\beta}\boldsymbol{X}_{i j}+\boldsymbol{\Lambda}\boldsymbol{\omega}_{i j} +\int_{\mathcal{T}} \boldsymbol{C}(t)\boldsymbol{Z}_{i j}(t)d t  +\boldsymbol{\epsilon}_{i j},
\end{aligned}
\end{equation}
where $\boldsymbol{\beta}$ is a $d_1 \times (d_2+1)$ regression coefficient matrix  including intercept, $\boldsymbol{\Lambda}$ is a $d_1 \times r_3$ loading matrix,   $\boldsymbol{C}(t)$ are coefficient functions which are defined on continuous domains,  modeling functional covariate effects. $\boldsymbol{\epsilon}_{i j}$ is a $d_1 \times 1$ random vector of error measurements with distribution $N(\mathbf{0}, \boldsymbol{\Psi})$, where $\boldsymbol{\Psi}=\operatorname{diag}(\psi_1, \ldots, \psi_{d_1})$ is a $d_1 \times d_1$ diagonal matrix, and $\boldsymbol{\epsilon_{i j}}$ is assumed independent with $\boldsymbol{\omega_{i j}}$, $\boldsymbol{X}_{i j}$ and $\boldsymbol{Z}_{i j}(t)$.

In Eq. \eqref{eq: PFSEM}, both the functional covariates $\boldsymbol{Z}_{ij}(t)$ and coefficient functions $\boldsymbol{C}(t)$ are defined on continuous domains, representing infinite-dimensional objects. To transform these into computationally tractable finite-dimensional problems, we employ basis expansion techniques for dimension reduction. Let $\{b_{k}(t)\}_{k=1}^{K_n}$ be a set of basis functions spanning a $K_n$-dimensional function space.  Any sufficiently smooth function can be approximated by a linear combination of these basis functions.
The parameter reduced form of Eq. \eqref{eq: PFSEM} then becomes:
\begin{equation}\label{eq: PFSEM_reduced}
\begin{aligned}
Y_{i j p}=\sum_{q=0}^{d_2}\beta_{p q}{X}_{i j q}+\sum_{r=1}^{r_3}\lambda_{p r}\omega_{i j r} +\sum_{h=1}^{d_3}\sum_{k=1}^{K_n} \phi_{p h k} Z_{i j h k} +\epsilon_{i j p},
\end{aligned}
\end{equation}
where $Z_{ijhk}=\int b_{k}(t) Z_{i j h}(t) d t$  are elements of the known design matrix, $\phi_{p h k}=\int b_{k}(t) C_{p h}(t) d t$ are finite-dimensional parameters to be estimated. 

The basis expansion technique approximates elements in the infinite-dimensional function space as elements in a finite-dimensional subspace spanned by basis functions, where $Z_{ijhk}$ represents the projection of functional data $Z_{ijh}(t)$ onto the basis function direction, and $\phi_{phk}$ represents the projection of the coefficient function $C_{ph}(t)$ onto the same direction. This approach provides a regularization method by controlling the number of basis functions $K_n$ to balance model complexity and goodness-of-fit. Additionally, it transforms integration operations into finite summations, enabling estimation with standard statistical software. Ultimately, this dimension reduction successfully converts the infinite-dimensional functional regression problem into a finite-dimensional multiple regression problem, significantly reducing computational complexity while preserving the essential characteristics of the functional data. The mathematical derivation is provided in Appendix \ref{sec:FDR_SI}.

We assume $\boldsymbol{\eta}_{i j}$ and $\boldsymbol{\xi}_{i j}$ satisfies the following nonlinear structural equation: 
\begin{equation}\label{eq: structural_equation}
\boldsymbol{\eta}_{i j}=\boldsymbol{\Gamma}_i \boldsymbol{G}\left(\boldsymbol{\xi}_{i j}\right)+\boldsymbol{\delta}_{i j},
\end{equation}
where $\boldsymbol{G}(\boldsymbol{\xi}_{i j})=(g_1(\boldsymbol{\xi}_{i j}), \ldots, g_b(\boldsymbol{\xi}_{i j}))^\top$ is a $b \times 1$ nonzero vector-valued function with differentiable functions $g_1,\cdots, g_b$ ($b\geq r_2$); $\boldsymbol{\Gamma}_i(r_1 \times r_2)$ is structural coefficient matrix denoting the effect of  $\boldsymbol{\xi}_{i j}$ on $ \boldsymbol{\eta}_{i j}$; $\boldsymbol{\delta}_{i j}(r_1 \times 1)$  is a residual vector independent of $\boldsymbol{\xi}_{ij}$. 

Modeling structural heterogeneity via the following equation: 
\begin{equation}\label{eq: structural_heterogeneity}	
\boldsymbol{\Gamma}_i=\boldsymbol{O}_i\boldsymbol{B}+\boldsymbol{U}_i\boldsymbol{v}_i,
\end{equation}
where $\boldsymbol{O}_i$ is an $r_1 \times 1$ vector of individual-level covariates useful for explaining $\boldsymbol{\Gamma}_i, \boldsymbol{B}$ is a $1 \times r_2$ vector of regression coefficients, $\boldsymbol{U}_i$ is an $r_1 \times 1$ indicator vector of $0$s and $1$s, and $\boldsymbol{v}_i=(v_{i 1}, \ldots, v_{i r_2})$ is a $1 \times r_2$ random vector independent of $\boldsymbol{\delta}_{i j}$ and $\boldsymbol{\xi}_{i j}$, we assume $\boldsymbol{v}_i\sim N(\mathbf{0}, \mathbf{\Upsilon})$, where $\mathbf{\Upsilon}$ is an $r_2 \times r_2$ covariance matrix representing the covariation of structural parameters due to  unobserved individual-level variables. 

In conventional SEMs, it is commonly assumed that $\boldsymbol{\xi}_{i j}$ is independent of $\boldsymbol{\xi}_{i l}$, and $\boldsymbol{\delta}_{i j}$ is also independent of $\boldsymbol{\delta}_{i l}$ for $j \neq l$.  Unlike conventional SEM formulations, our framework explicitly incorporates both between-variable and adjacent-time covariance structures for $\boldsymbol{\xi}_{ij}$ and $\boldsymbol{\delta}_{ij}$. We model these dependencies through the Linear Model of Coregionalization (LMC) \citep{song2012bayesian}:
$$
\boldsymbol{\xi}_{i j}=\boldsymbol{A}_1 \mathbf{w}_{i j}, \quad \boldsymbol{\delta}_{i j}=\boldsymbol{A}_2 \mathbf{v}_{i j},
$$
where $\boldsymbol{A}_1$ and $\boldsymbol{A}_2$ are $r_2 \times r_2$ and $r_1 \times r_1$ upper triangular matrices, respectively; $\mathbf{w}_{i j}$ and $\mathbf{v}_{i j}$ are independent zero-mean and unit-variance random vectors of dimensions $r_2$ and $r_1$, respectively. 
Let $\boldsymbol{C}_i(\boldsymbol{\alpha}_{i 1})$ and $\boldsymbol{D}_i(\boldsymbol{\alpha}_{i 2})$ be the $J \times J$ adjacent time covariance matrices of vectors
$\tilde{\mathbf{w}}_{i k}=(\mathrm{w}_{i 1 k}, \ldots, \mathrm{w}_{i J k})^\top$ and $\tilde{\mathbf{v}}_{i m}=(\mathrm{v}_{i 1 m}, \ldots, \mathrm{v}_{i J m})^\top$, respectively,  in which $\mathrm{w}_{i j k}$ and $\mathrm{v}_{i j m}$ denote the $k$ th and $m$ th elements of  $\mathbf{w}_{i j}$ and $\mathbf{v}_{i j}$ ($k=1, \ldots, r_2$, $m=1, \ldots, r_1$); $\boldsymbol{\alpha}_{i 1}$ and $\boldsymbol{\alpha}_{i 2}$ are the corresponding vectors of parameters in the covariance structures.  Let $\boldsymbol{\xi}_i=(\boldsymbol{\xi}_{i 1}^\top, \ldots, \boldsymbol{\xi}_{i J}^\top)^\top$ and $\boldsymbol{\delta}_i=(\boldsymbol{\delta}_{i 1}^\top, \ldots, \boldsymbol{\delta}_{i J}^\top)^\top$.  Then, the covariance matrices $\boldsymbol{\Sigma}_{i \xi}$ and $\boldsymbol{\Sigma}_{i \delta}$ of $\boldsymbol{\xi}_i$ and $\boldsymbol{\delta}_i$ can be expressed as
\begin{equation}\label{eq: covariance_matrices}
\boldsymbol{\Sigma}_{i \xi}=\boldsymbol{C}_i\left(\alpha_{i 1}\right) \otimes \boldsymbol{T}_{\xi}, \quad \boldsymbol{\Sigma}_{i \delta}=\boldsymbol{D}_i\left(\alpha_{i 2}\right) \otimes \boldsymbol{T}_\delta,
\end{equation}
where  $\boldsymbol{T}_{\xi}=\boldsymbol{A}_1 \boldsymbol{A}_1^\top$ and $\boldsymbol{T}_\delta=\boldsymbol{A}_2 \boldsymbol{A}_2^\top$ represent the between-variable covariances of $\boldsymbol{\xi}_{i j}$ and $\boldsymbol{\delta}_{i j}$, respectively.  The adjacent time dependency structure is formulated using the conditional autoregressive (CAR) model \citep{Besag1974SpatialIA}, specified as
\begin{equation}\label{eq: CAR}
\boldsymbol{C}_i\left(\boldsymbol{\alpha}_{i 1}\right)=\tau_{i \xi}^2\left(\boldsymbol{I}_{J}-\rho_{i \xi} \boldsymbol{H}_{i \xi}\right)^{-1}, \quad \boldsymbol{D}_i\left(\boldsymbol{\alpha}_{i 2}\right)=\tau_{i \delta}^2\left(\boldsymbol{I}_{J}-\rho_{i \delta} \boldsymbol{H}_{i \delta}\right)^{-1},
\end{equation}
where $\rho_{i \xi}$ and $\rho_{i \delta}$ quantify temporal dependence between adjacent time points, $\tau_{i \xi}^2$ and  $\tau_{i \delta}^2$ are parameters proportional to the conditional variances of the elements in $\tilde{\mathbf{w}}_{i k}$ and $\tilde{\mathbf{v}}_{i m}$ respectively,  given their neighbors; $\boldsymbol{H}_{i \xi}$ and $\boldsymbol{H}_{i \delta}$ are $J \times J$ matrices containing the neighboring information,  where element $h_{j l}=1$ indicates that time $l$ is adjacent to time $j$ and $h_{j l}=0$ otherwise. This study considers exclusively first-order temporal adjacency ($l = j \pm 1$), with $\boldsymbol{H}_{i\xi} = \boldsymbol{H}_{i\delta} = (h_{jl})_{J \times J}$. In CAR model, $\boldsymbol{\alpha}_{i 1}=(\tau_{i \xi}^2, \rho_{i \xi})^\top$ and $\boldsymbol{\alpha}_{i 2}=(\tau_{i \delta}^2, \rho_{i \delta})^\top$.

We proposed the partially functional dynamic structural equation model (PFDSEM) that offers several key methodological advantages for analyzing complex longitudinal data with functional and scalar covariates, as well as latent structures. A major strength of our approach lies in its ability to simultaneously incorporate scalar, functional, and latent variables within a unified modeling framework.

First, the functional data processing approach via basis expansion offers significant advantages for handling functional covariates within statistical modeling frameworks. Its primary strength lies in transforming infinite-dimensional problems into tractable finite-dimensional forms by projecting both functional covariates and coefficient functions onto a finite basis system \citep{Ramsay2005FDA} :  
i) This approach provides substantial computational benefits while maintaining modeling flexibility through the adaptability of basis functions to different data characteristics, Fourier bases for periodic patterns \citep{Ramsay2005FDA}, B-splines for smooth non-periodic relationships \citep{Hastie2009ESL}, and wavelets for capturing local features \citep{Hastie2009ESL}; 
ii) Unlike discretization or grid analysis that treats time points as separate covariates and fails to capture functional relationships, basis expansion preserves the continuous nature of the data while reducing dimensionality \citep{Ramsay2005FDA}; 
iii) In contrast to two-step approaches that perform pre-smoothing independently of outcome modeling \citep{Wood2017GAM}, our method simultaneously handles functional representation and relationship estimation within a unified framework, thereby reducing potential bias; 
iv) Differing from Functional Principal Component Analysis (FPCA), which prioritizes explaining covariate variance \citep{Hastie2009ESL}, basis expansions directly target covariate-response relationships, typically yielding more parsimonious and predictive models \citep{Wood2017GAM}; 
v) The method ensures clear interpretability through reconstructible functional coefficients and allows seamless integration into broader multivariate modeling frameworks. This makes it particularly valuable for complex models incorporating functional, latent, and observed variables.    

Furthermore, the PFDSEM introduces a dynamic structure that melds a conditional autoregressive (CAR) specification, a type of Markov random field prior, with a linear model of coregionalization (LMC). This combination provides a statistically principled mechanism for modeling complex dependence: the CAR structure captures the inter-temporal smoothing and evolution of the process, while the LMC flexibly models the inter-variable covariances, generating a rich, time-varying correlation structure for the latent vector. This explicit separation of dependency types enhances model interpretability and provides inherent regularization. The PFDSEM is further extended to include unit-level random effects, accounting for heterogeneity in structural parameters, and to support nonlinear latent relationships. 


\subsection{Model Identifiability}\label{sec:MIC} 

To ensure model identifiability, we need the following essential constraints.

\begin{condition}(Factor Loading Matrix and Covariance Matrix  Constraints)
\label{con:1}
The first non-zero element of each column in matrix $\boldsymbol{\Lambda}$ is set to 1 to avoid rotational indeterminacy\citep{Joreskog_1969}, and fix $\tau_{i\xi}^2 = \tau_{i\delta}^2 = 1$ to avoid scale indeterminacy\citep{Besag1974SpatialIA}.
\end{condition}

\begin{condition}
(CAR Model Parameter Constraints)
\label{con:2}
 $\rho_{i \xi}$ and $\rho_{i \delta}$ must satisfy the following constrained conditions $\rho_{i \xi}^L<\rho_{i \xi}<\rho_{i \xi}^U$ and $\rho_{i \delta}^L<\rho_{i \delta}<\rho_{i \delta}^U$ to ensure that $\boldsymbol{C}_i$ and $\boldsymbol{D}_i$ are positive definite, where $\rho_{i \xi}^L=\min (\frac{1}{\tau_{\xi T_i}}, 0), \rho_{i \xi}^U=\frac{1}{\tau_{\xi 1}}, \rho_{i \delta}^L=\min (\frac{1}{\tau_{\delta \tau_i}}, 0), \rho_{i \delta}^U=\frac{1}{\tau_{\delta 1}}$, and $\tau_{\xi 1} \geq \cdots \geq \tau_{\xi J}$ and $\tau_{\delta 1} \geq \cdots \geq \tau_{\delta J}$ are the eigenvalues of neighborhood matrices $\boldsymbol{H}_{i \xi}$ and $\boldsymbol{H}_{i \delta}$, respectively \citep{song2012bayesian}.
\end{condition}

\begin{condition}
(Linear Model of Coregionalization Identification)
\label{con:3}
For the LMC parameters $\boldsymbol{A}_1$ and $\boldsymbol{A}_2$ to be identified, we require: 
\begin{enumerate}
    \item[(i)] The temporal covariance matrices $\boldsymbol{C}_i(\boldsymbol{\alpha}_{i1})$ and $\boldsymbol{D}_i(\boldsymbol{\alpha}_{i2})$ are identified under Condition \ref{con:2} together with the scale constraints $\tau_{i\xi}^2=\tau_{i\delta}^2=1$ in Condition \ref{con:1}.
    \item[(ii)] $\boldsymbol{A}_1$ and $\boldsymbol{A}_2$ are specified as upper triangular matrices with positive diagonal elements (i.e., Cholesky decompositions of $\boldsymbol{T}_{\xi}$ and $\boldsymbol{T}_{\delta}$, respectively).
\end{enumerate}
\end{condition}
 
Condition \ref{con:1} addresses two classical sources of indeterminacy in latent variable models. Fixing the first nonzero element of each column of the factor loading matrix $\boldsymbol{\Lambda}$ to 1 is a standard practice in structural equation modeling \citep{Joreskog_1969}. This constraint eliminates rotational indeterminacy by anchoring the scale of each latent variable to a specific manifest indicator, thereby allowing unique estimation of the loadings. Furthermore, setting $\tau_{i\xi}^2 = \tau_{i\delta}^2 = 1$ resolves scale indeterminacy in the conditional autoregressive (CAR) components \citep{Besag1974SpatialIA}. This restriction fixes the conditional variance parameters to a constant, preventing arbitrary scaling of the adjacent-time covariance matrices $\boldsymbol{C}_i$ and $\boldsymbol{D}_i$, while still allowing the dependence strength to be flexibly captured through the parameters $\rho_{i\xi}$ and $\rho_{i\delta}$. 

Condition \ref{con:2} ensures the positive definiteness of the covariance matrices $\boldsymbol{C}_i$ and $\boldsymbol{D}_i$, which is a necessary condition for valid covariance modeling and numerical stability during estimation. The constraints on $\rho_{i\xi}$ and $\rho_{i\delta}$ are derived from the spectral properties of the neighborhood matrices $\boldsymbol{H}_{i\xi}$ and $\boldsymbol{H}_{i\delta}$ \citep{song2012bayesian}. Specifically, the bounds $\rho_{i\xi}^L < \rho_{i\xi} < \rho_{i\xi}^U$ and $\rho_{i\delta}^L < \rho_{i\delta} < \rho_{i\delta}^U$ guarantee that the precision matrices $(\boldsymbol{I} - \rho_{i\xi} \boldsymbol{H}_{i\xi})$ and $(\boldsymbol{I} - \rho_{i\delta} \boldsymbol{H}_{i\delta})$ remain invertible and that the resulting covariance structures are positive definite. This ensures that the temporal dependence parameters yield a valid probabilistic structure while still accommodating a wide range of dependence scenarios from negative to positive spatial autocorrelation. 

Condition \ref{con:3} establishes that the parameters in the Linear Model of Coregionalization are uniquely determined. The Kronecker product structure $\boldsymbol{\Sigma}_{i\xi} = \boldsymbol{C}_i(\boldsymbol{\alpha}_{i1}) \otimes \boldsymbol{T}_{\xi}$ separates temporal dynamics (captured by $\boldsymbol{C}_i$) from cross-sectional dependence (captured by $\boldsymbol{T}_{\xi}$). Part (i) guarantees that $\boldsymbol{C}_i$ is already identified from the CAR specification under the previous constraints. Given $\boldsymbol{C}_i$, the cross-sectional covariance $\boldsymbol{T}_{\xi} = \boldsymbol{A}_1\boldsymbol{A}_1^\top$ would normally be identified only up to an orthogonal rotation $\boldsymbol{A}_1\boldsymbol{Q}$ ($\boldsymbol{Q}\boldsymbol{Q}^\top=\boldsymbol{I}$). Part (ii) eliminates this rotational indeterminacy by requiring $\boldsymbol{A}_1$ to be upper triangular with positive diagonal entries, the unique Cholesky factor of $\boldsymbol{T}_{\xi}$ \citep{Pourahmadi1999}. The same reasoning applies to $\boldsymbol{A}_2$ and $\boldsymbol{T}_{\delta}$. Hence, under Condition \ref{con:3}, all LMC parameters become fully identifiable.

\section{Bayesian Inference }\label{sec:BE} 

Statistical inference for the proposed PFDSEM is conducted under the Bayesian paradigm. This framework provides a coherent probabilistic approach for estimating complex hierarchical models, where all unknown quantities, including parameters and latent variables, are treated as random variables. Uncertainty regarding these unknowns is formally quantified through their probability distributions \citep{gelman2013bayesian}.

Let $\boldsymbol{\theta}_1$ contain all unknown distinct parameters in $\boldsymbol{\beta}, \boldsymbol{\Lambda}$, $\boldsymbol{\Phi}$ and $\boldsymbol{\Psi}$  associated with Eq. \eqref{eq: PFSEM_reduced}, $\boldsymbol{\theta}_2$ contain all unknown distinct parameters in $\boldsymbol{B}, \boldsymbol{\Upsilon}, \boldsymbol{A}_1, \boldsymbol{A}_2, \boldsymbol{\alpha}_1, \boldsymbol{\alpha}_2$ associated with Eqs. \eqref{eq: structural_equation}-\eqref{eq: CAR}, and $\boldsymbol{\theta}=\left\{\boldsymbol{\theta}_1, \boldsymbol{\theta}_2\right\}$, where $\boldsymbol{B}=$ $(B_{1}, \ldots, B_{r_2}), \boldsymbol{\alpha}_1= ( \boldsymbol{\alpha}_{11}, \ldots, \boldsymbol{ \alpha}_{n1} )$, and $\boldsymbol{\alpha}_2 = ( \boldsymbol{\alpha}_{12}, \ldots, \boldsymbol{\alpha}_{n2} )$. Let $\boldsymbol{Y} = \left\{ \boldsymbol{Y}_{i j}: i=1, \ldots, n, j=1, \ldots, J \right\}$, $ \boldsymbol{F}_{i(1)} = (\bm{\eta}_{i 1}, \ldots, \bm{\eta}_{i J})$, $\boldsymbol{F}_{i(2)} = ( \boldsymbol{\xi}_{i 1}, \ldots, \boldsymbol{\xi}_{i J} )$, $\boldsymbol{F}_1 = ( \boldsymbol{F}_{1(1)}, \ldots, \boldsymbol{F}_{n(1) } )$, $\boldsymbol{F}_2 = (\boldsymbol{F}_{1(2)}, \ldots, \boldsymbol{F}_{n(2) })$, $\boldsymbol{F} =(\boldsymbol{F}_1^\top, \boldsymbol{F}_2^\top )^\top$, $\boldsymbol{X} = \left\{\boldsymbol{x}_{i j}: i=1, \ldots, n, j=1, \ldots, J\right\}$, $\boldsymbol{O}= \left\{\bm{O}_i: i=1, \ldots, n\right\}$ and $\ \ \boldsymbol{U} = \left\{  \boldsymbol{U}_i:  i=1, \ldots, n\right\}$.

The foundation of Bayesian inference is Bayes' theorem, which updates prior beliefs about the parameters, encoded in the joint prior density \(p(\boldsymbol{\theta}, \boldsymbol{\Gamma})\), with the information contained in the observed data via the likelihood function \(p(\boldsymbol{D} \mid \boldsymbol{\theta}, \boldsymbol{\Gamma})\). This yields the posterior distribution:
\begin{equation}\label{eq:bayes_thm}
 p(\boldsymbol{\theta}, \boldsymbol{\Gamma} \mid \boldsymbol{D}) = \frac{p(\boldsymbol{D} \mid \boldsymbol{\theta}, \boldsymbol{\Gamma}) \, p(\boldsymbol{\theta}, \boldsymbol{\Gamma})}{\int p(\boldsymbol{D} \mid \boldsymbol{\theta}, \boldsymbol{\Gamma}) \, p(\boldsymbol{\theta}, \boldsymbol{\Gamma}) \, d\boldsymbol{\theta} d\boldsymbol{\Gamma}} \propto p(\boldsymbol{D} \mid \boldsymbol{\theta}, \boldsymbol{\Gamma}) \, p(\boldsymbol{\theta}, \boldsymbol{\Gamma}).
\end{equation}
The posterior distribution \(p(\boldsymbol{\theta}, \boldsymbol{\Gamma} \mid \boldsymbol{D})\) represents the complete solution to the inference problem, synthesizing all available information.

The PFDSEM is a complex hierarchical model, direct evaluation of the marginal posterior \(p(\boldsymbol{\theta}, \boldsymbol{\Gamma} \mid \boldsymbol{D})\) is analytically intractable due to the high-dimensional integration over the latent space. To facilitate computation, we adopt a data augmentation approach \citep{Tanner01061987}, treating the latent variables \(\boldsymbol{F}\) as auxiliary parameters. This leads to the augmented joint posterior distribution:
\begin{equation}\label{eq:augmented_posterior}
 p(\boldsymbol{\theta}, \boldsymbol{\Gamma}, \boldsymbol{F} \mid \boldsymbol{D}) \propto p(\boldsymbol{D} \mid \boldsymbol{\theta}, \boldsymbol{\Gamma}, \boldsymbol{F}) \, p(\boldsymbol{F} \mid \boldsymbol{\theta}, \boldsymbol{\Gamma}) \, p(\boldsymbol{\theta}, \boldsymbol{\Gamma}),
\end{equation}
where the likelihood \(p(\boldsymbol{D} \mid \boldsymbol{\theta}, \boldsymbol{\Gamma}, \boldsymbol{F})\) is now conditional on the latent states. This formulation defines a hierarchical model structure that is more amenable to Gibbs sampling. We employ a Markov Chain Monte Carlo (MCMC) method, specifically a Gibbs sampler \citep{Geman1984Gibbs, robert2004monte}, to draw samples from the joint posterior distribution in Eq.\eqref{eq:augmented_posterior}.  A Gibbs sampler \citep{Geman1984Gibbs} generates sequences from $p(\boldsymbol{\theta}, \boldsymbol{\Gamma}, \mathbf{F} \mid \mathbf{D})$ using initial values $(\boldsymbol{\theta}^{(0)}, \boldsymbol{\Gamma}^{(0)},\mathbf{F}^{(0)})$ through the following iterations:  i) Sample $\mathbf{F}^{(k+1)} \sim p(\mathbf{F} \mid \boldsymbol{\theta}^{(k)}, \boldsymbol{\Gamma}^{(k)}, \mathbf{D})$;  ii) Sample $\boldsymbol{\theta}^{(k+1)} \sim p(\boldsymbol{\theta} \mid \mathbf{F}^{(k+1)}, \boldsymbol{\Gamma}^{(k)}, \mathbf{D})$;  iii) Sample $\boldsymbol{\Gamma}^{(k+1)} \sim p(\boldsymbol{\Gamma} \mid \mathbf{F}^{(k+1)}, \boldsymbol{\theta}^{(k+1)}, \mathbf{D})$. Under regularity conditions \citep{Geman1984Gibbs}, $(\boldsymbol{\theta}^{(K)}, \boldsymbol{\Gamma}^{(K)},\mathbf{F}^{(K)})$ converges in distribution to the joint posterior for sufficiently large $K$. Posterior distributions are approximated using $\{( \boldsymbol{\theta}^{(k)}, \boldsymbol{\Gamma}^{(k)},\mathbf{F}^{(k)}): k=K+1, \ldots, K+T\}$ with large $T$.

Reliable inference requires that the MCMC chains have converged to the target distribution. 
We assess convergence using a multi-chain diagnostic procedure \citep{gelman2013bayesian}. 
We run multiple independent chains (typically three) from dispersed starting points and monitor the ``Estimated Potential Scale Reduction (EPSR)'' values \citep{Brooks01121998}.  
Convergence is achieved when all the EPSR values fall below $1.2$.  After convergence, a large number of observations simulated from the joint posterior distribution $p(\boldsymbol{\theta},  \boldsymbol{\Gamma},\boldsymbol{F} | \boldsymbol{D})$ can be used to produce the bayesian estimates and the standard errors estimates of the unknown parameters and latent variables with the help of the corresponding posterior means and their posterior covariance matrices.  This Bayesian approach offers several advantages for the PFDSEM, it naturally propagates uncertainty from all model levels (measurement, structural, spatial-temporal), avoids reliance on asymptotic approximations, and provides a direct probabilistic interpretation for all estimates \citep{Lee2007SEM}.

\subsection{Priors}\label{sec:PS} 

The specification of prior distributions $p(\boldsymbol{\theta}, \boldsymbol{\Gamma})$ is a critical component of Bayesian analysis. Well-chosen priors can regularize estimation, incorporate domain knowledge, and ensure proper posterior distributions. For the PFDSEM, we adopt conditionally conjugate priors where possible to facilitate Gibbs sampling, and we use weakly informative or reference priors when substantive prior information is limited \citep{gelman2006data}.

For the regression coefficients $\boldsymbol{\beta}= (\boldsymbol{\beta}_{0}^\top,\cdots, \boldsymbol{\beta}_{p}^\top, \cdots,\boldsymbol{\beta}_{d_1}^\top)^\top_{d_1 \times (d_2+1)}$, factor loadings $\boldsymbol{\Lambda}= (\boldsymbol{\Lambda_{1}}, \cdots, \boldsymbol{\Lambda_{p}},$ $ \cdots,\boldsymbol{\Lambda}_{d_1} )_{r_3 \times d_1}$, and residual variances $\boldsymbol{\Psi}$ in the measurement equation (Eq. \eqref{eq: PFSEM_reduced}), we specify conditionally conjugate priors \citep{Lee2007SEM}:
\begin{equation}\label{eq: Parameters_Priors}
\begin{aligned}
p(\boldsymbol{\beta}_{p}) \stackrel{D}{=} N_{d_2+1}\left(\boldsymbol{\mu}_0, \boldsymbol{\Sigma}_{0 \beta_{p}}\right), \quad
p\left(\boldsymbol{\Lambda}_p \mid \psi_p\right) \stackrel{D}{=} N_{r_3}\left(\boldsymbol{\Lambda}_{0 p}, \psi_p \boldsymbol{\Sigma}_{0 \Lambda_{p}}\right), \quad
 p\left(\psi_p^{-1}\right) \stackrel{D}{=} \operatorname{Gamma}\left(\alpha_{0 \psi p}, \beta_{0 \psi p}\right),
 \end{aligned}
\end{equation}
where $\boldsymbol{\mu}_0, \boldsymbol{\Sigma}_{0 \beta_{p}}, \boldsymbol{\Lambda}_{0 p}, \boldsymbol{\Sigma}_{0 \Lambda p}, \alpha_{0 \psi p}$, and $\beta_{0 \psi p}$ are prespecified hyperparameters. These choices yield standard full conditional distributions: normal for $\boldsymbol{\beta}_p$ and $\boldsymbol{\Lambda}_p$, and Gamma for $\psi_p^{-1}$, facilitating straightforward Gibbs sampling.  
The conditional prior for $\boldsymbol{\Lambda}_p$ given $\psi_p$ is standard in Bayesian factor analysis and SEM \citep{Lee2007SEM}, it ensures that the scale of the loadings is appropriately linked to the residual variance, and together with the gamma prior on $\psi_p^{-1}$. 
We assume independence between $(\psi_p, \boldsymbol{\Lambda}_p)$ and $(\psi_m, \boldsymbol{\Lambda}_m)$ for $p \neq m$. When prior information is limited, $\boldsymbol{\Sigma}_{0 \beta_{p}}$ and $\boldsymbol{\Sigma}_{0 \Lambda p}$ are set to large values to create weakly informative priors, while 
the gamma hyperparameters are chosen to be small ($\alpha_{0\psi p}$ and  $\beta_{0\psi p}$) to provide weak information.

For the functional coefficients $\boldsymbol{\Phi}_p=(\boldsymbol{\phi}_{p 1},\cdots,\boldsymbol{\phi}_{p h},\cdots,\boldsymbol{\phi}_{p d_3})^\top_{d_3 \times K_n}$ in Eq. \eqref{eq: PFSEM_reduced}, we employ Bayesian P-splines \citep{Lang2004BayesianPSplines}. This approach provides a flexible nonparametric representation while controlling overfitting through a smoothness penalty encoded in the prior:
\begin{equation}\label{eq: Parameters_Priors2}
\boldsymbol{\phi}_{p h} \sim \operatorname{N_{K_n}}(\mathbf{0}, \theta_{p h}^2 \boldsymbol{U}^{-1}), 
\end{equation}
where $\boldsymbol{U}$ is the second-order difference penalty matrix that penalizes roughness in the coefficient function. 
This prior is a multivariate normal with precision matrix proportional to $\boldsymbol{U}$, which penalizes roughness in the coefficient function; it is equivalent to assuming a random walk prior on the B-spline coefficients and is widely used in semiparametric regression \citep{Lang2004BayesianPSplines, Wood2017GAM}. 
Following \cite{Klein2016ScaleDependentPriors}, we assign a hyperprior to the smoothing parameter:
$$
p\left(\theta_{p h}^{-2}\right) \stackrel{D}{=} \operatorname{Gamma}\left(\alpha_{0 \theta_{p h}}, \beta_{0 \theta_{p h}}\right), 
$$
where $\alpha_{0 \theta_{p h}}$ and $\beta_{0 \theta_{p h}}$ are prespecified. This scale-dependent prior provides adaptive smoothing, allowing the data to determine the appropriate level of smoothness for each coefficient function.

For the structural heterogeneity parameters $\boldsymbol{B}$ (fixed effects)  and $\boldsymbol{\Upsilon}$ (random effects covariance)  in Eq. \eqref{eq: structural_heterogeneity},
we specify conjugate priors:
\begin{equation}\label{eq: Parameters_Priors3}
p(\boldsymbol{B}) \stackrel{D}{=} N_{r_2}\left(\boldsymbol{B}_0, \boldsymbol{\Sigma}_{0 B}\right),\quad p\left(\mathbf{\Upsilon}^{-1}\right) \stackrel{D}{=} \text { Wishart }_{r_2}\left(\boldsymbol{\mu}_{0 \Upsilon}, \boldsymbol{R}_{0 \Upsilon}\right),
\end{equation}
where $\boldsymbol{B}_0, \boldsymbol{\Sigma}_{0 B}, \boldsymbol{\mu}_{0 \Upsilon}$, and $\boldsymbol{R}_{0 \Upsilon}$ are prespecified. 
The normal prior for $\boldsymbol{B}$ is conditionally conjugate, and the Wishart prior for $\boldsymbol{\Upsilon}^{-1}$ is the standard conjugate choice for covariance matrices in hierarchical models \citep{gelman2013bayesian}.  If $\boldsymbol{\Upsilon}$ is assumed diagonal with elements $\Upsilon_1, \ldots, \Upsilon_{r_2}$ (e.g., for independent random effects) ,
we use independent inverse Gamma priors: $p(\Upsilon_l^{-1}) \stackrel{D}{=} \operatorname{Gamma}(\alpha_{0 \Upsilon}, \beta_{0 \Upsilon})$ for $l=1, \ldots, r_2$, which maintains conjugacy while allowing separate variance components. where $\alpha_{0 \Upsilon}$ and $\beta_{0 \Upsilon}$ are given hyperparameters.

For the between-variable covariance matrices $\boldsymbol{T}_{\xi}$ and $\boldsymbol{T}_\delta$ (Eq. \eqref{eq: covariance_matrices}), we assign conjugate Wishart priors:
\begin{equation}\label{eq: Parameters_Priors4}
\begin{gathered}
p\left(\boldsymbol{T}_{\xi}^{-1}\right) \stackrel{D}{=} \text { Wishart }_{r_2}\left(\boldsymbol{\mu}_{0 \xi}, \boldsymbol{R}_{0 \xi}\right),  \quad p\left(\boldsymbol{T}_\delta^{-1}\right) \stackrel{D}{=} \text { Wishart }_{r_1}\left(\boldsymbol{\mu}_{0 \delta}, \boldsymbol{R}_{0 \delta}\right),
\end{gathered}
\end{equation}
where $\boldsymbol{\mu}_{0 \xi}, \boldsymbol{\mu}_{0 \delta}, \boldsymbol{R}_{0 \xi}$, and $\boldsymbol{R}_{0 \delta}$ are prespecified hyperparameters.
The Wishart distribution is the natural conjugate prior for the precision matrix of a multivariate normal distribution, leading to Wishart full conditionals that facilitate Gibbs sampling \citep{Lee2007SEM}. Hyperparameters are chosen to be weakly informative (e.g., degrees of freedom $\mu_{0\xi}=r_2+5$, scale matrix $\boldsymbol{R}_{0\xi}$ diagonal with moderate values).

Following Condition \ref{con:2}, we specify uniform priors over their valid ranges for $\rho_{i \xi}$ and $\rho_{i \delta}$ in Eq. \eqref{eq: CAR}:
\begin{equation}\label{eq: Parameters_Priors5}
p\left(\rho_{i \xi}\right) \stackrel{D}{=} \operatorname{Uniform}\left(\rho_{i \xi}^L, \rho_{i \xi}^U\right), \quad p\left(\rho_{i \delta}\right) \stackrel{D}{=} \operatorname{Uniform}\left(\rho_{i \delta}^L, \rho_{i \delta}^U\right) .
\end{equation}
These proper but non-informative priors respect the constraints required for positive definiteness of the CAR covariance matrices \citep{Besag1974SpatialIA}. 
They also reflect the lack of prior knowledge about the strength of temporal dependence beyond the admissible range.

All hyperparameters are either set to default weakly informative values or chosen based on empirical Bayes considerations in the real data analysis. Prior sensitivity is assessed in Section \ref{sec:SS} by varying hyperparameter choices and examining their impact on posterior estimates.

\subsection{Conditional Distributions}\label{sec:CD} 

To implement the Gibbs sampler, we require the full conditional distributions of all unknown quantities. These distributions are derived from the joint posterior and, where possible, take standard forms that permit direct sampling. For parameters whose conditionals are non-standard, we employ Metropolis-Hastings steps within the Gibbs iteration. Detailed derivations are provided in Appendix \ref{sec:CD_SI}, here we briefly outline the key results.

From the definitions of $\mathbf{F}$ and $\mathbf{D}$, the conditional distribution:
\begin{equation}\label{eq: CD_F}
\begin{aligned}
&p(\boldsymbol{F} \mid \boldsymbol{Y},\boldsymbol{\Gamma}, \boldsymbol{A}_2, \boldsymbol{A}_1, \boldsymbol{\alpha}_1, \boldsymbol{\alpha}_2)\propto p\left(\boldsymbol{Y} \mid \boldsymbol{F}, \boldsymbol{\beta}, \boldsymbol{\Lambda}, \boldsymbol{\Psi},\boldsymbol{\Phi},\boldsymbol{X},\boldsymbol{Z}\right) p\left(\boldsymbol{F}_1 \mid \boldsymbol{F}_2,  \boldsymbol{\Gamma}, \boldsymbol{A}_2, \boldsymbol{\alpha}_2\right) p\left(\boldsymbol{F}_2 \mid \boldsymbol{A}_1, \boldsymbol{\alpha}_1\right)\\
 \propto&\left(2 \pi\right)^{-\frac{N\left(d_1+r_3\right)}{2}} \left|\prod_{p=1}^{d_1} \psi_p\right|^{-\frac{N}{2}} \left|\prod_{i=1}^n \boldsymbol{C}_i\left(\boldsymbol{\alpha}_{i 1}\right)\right|^{-\frac{r_2}{2}}  \left|\prod_{i=1}^n \boldsymbol{D}_i\left(\boldsymbol{\alpha}_{i 2}\right)\right|^{-\frac{r_1}{2}} \left(\left|T_\xi\right| \left| T_\delta\right|\right)^{-\frac{N}{2}} \\
&\times \exp \Bigg\{-\frac{1}{2} \sum_{i=1}^n\Bigg(\sum_{j=1}^J \sum_{p=1}^{d_1} \psi^{-1}_p\left(y_{i j p}-\boldsymbol{\beta}_p \boldsymbol{X}_{i j}-\boldsymbol{\Lambda_p} \boldsymbol{\omega}_{i j}-\sum_{h=1}^{d_3} \boldsymbol{\phi}_{p h} \boldsymbol{Z}_{i j h}\right)^2\Bigg.\Bigg.\\
&\Bigg.\Bigg.+\operatorname{tr}\bigg(\boldsymbol{F}_{i(2)}^{\top} T_\xi^{-1} \boldsymbol{F}_{i(2)} \boldsymbol{C}_i\left(\boldsymbol{\alpha}_{i 1}\right)^{-1}+\boldsymbol{G}_i^{* T} \boldsymbol{\Xi}_i^{\top} T_\delta^{-1} \boldsymbol{\Xi}_i \boldsymbol{G}_i^* \boldsymbol{D}_i\left(\boldsymbol{\alpha}_{i 2}\right)^{-1}\bigg)\Bigg)\Bigg\}, 
\end{aligned}
\end{equation}
where $N=n \times J, $ $\boldsymbol{\Xi}_i=\left(\boldsymbol{I}_{r_2},-\boldsymbol{\Gamma}_i\right)$ and $\boldsymbol{G}_i^*=\left(\tilde{\boldsymbol{G}}_{i 1}, \ldots, \tilde{\boldsymbol{G}}_{i J}\right)$ with $\tilde{\boldsymbol{G}}_{i j}=(\boldsymbol{\eta}_{i j}^\top, \boldsymbol{G}(\boldsymbol{\xi}_{i j})^\top)^\top$. 

For each manifest variable $p=1,\ldots,d_1$, the conditional distributions of the regression coefficients and factor loadings are multivariate normal:
\begin{equation}\label{eq: CD_Beta-p}
p\left(\bm{\beta}_p \mid \bm{\Lambda},\bm{\Psi},\bm{\phi},\bm{F},\bm{Y},\bm{X},\bm{Z}\right) \stackrel{D}{=}  N_{d_2}\left(\bm{\mu}_{\beta_p}, \bm{\Sigma}_{\beta_p}\right),
\end{equation}
\begin{equation}\label{eq: CD_Lambda-p}
p\left(\bm{\Lambda}_p \mid \bm{\beta},\bm{\Psi},\bm{\phi},\bm{F},\bm{Y},\bm{X},\bm{Z}\right)\stackrel{D}{=}  N_{r_3}\left(\bm{\mu}_{\Lambda_p}, \bm{\Sigma}_{\Lambda_p}\right),
\end{equation}
where 
 $$
 \begin{aligned}     
&\bm{\Sigma}_{\beta_p}^{-1} = \bm{\Sigma}_{0\beta_p}^{-1} + \sum_{i=1}^{n}\sum_{j=1}^{J} \psi^{-1}_p \bm{X}_{ij}\bm{X}_{ij}^{\top},\quad  \bm{\Sigma}_{\Lambda_p}^{-1} =\psi^{-1}_p\bm{\Sigma}_{0\Lambda_p}^{-1} + \sum_{i=1}^{n}\sum_{j=1}^{J} \psi^{-1}_p \bm{\omega}_{ij}\bm{\omega}_{ij}^{\top},\\
&\bm{\mu}_{\beta_p} =\bm{\Sigma}_{\beta_p} \left( \bm{\Sigma}_{0\beta_p}^{-1}\bm{\mu}_0 + \sum_{i=1}^{n}\sum_{j=1}^{J} \psi^{-1}_p \left(y_{ijp} - \bm{\Lambda}_p\bm{\omega}_{ij} - \sum_{h=1}^{d_3}\bm{\phi}_{ph}\bm{Z}_{ijh}\right)\bm{X}_{ij} \right) ,\\
 &\bm{\mu}_{\Lambda_p} = \bm{\Sigma}_{\Lambda_p} \left( \psi^{-1}_p\bm{\Sigma}_{0\Lambda_p}^{-1}\bm{\Lambda}_{0p} + \sum_{i=1}^{n}\sum_{j=1}^{J} \psi^{-1}_p \left(y_{ijp} - \bm{\beta}_p\bm{X}_{ij} - \sum_{h=1}^{d_3}\bm{\phi}_{ph}\bm{Z}_{ijh}\right)\bm{\omega}_{ij} \right).
 \end{aligned}
 $$

The conditional distribution of $\boldsymbol{\Psi}$ combines the gamma prior with the likelihood contribution and the conditional prior for $\boldsymbol{\Lambda}$:
\begin{equation}\label{eq: CD_Psi}
\begin{aligned}
&p\left(\boldsymbol{\Psi} \mid \boldsymbol{\beta},  \boldsymbol{ \phi}, \boldsymbol{F},\boldsymbol{Y}, \boldsymbol{X}, \boldsymbol{Z}\right)\propto  p(\boldsymbol{\Lambda} \mid \boldsymbol{F}, \boldsymbol{Y},\boldsymbol{\Psi},\boldsymbol{ \phi}, \boldsymbol{\beta}) p\left(\boldsymbol{\Psi}\right) \\
\propto &\prod_{p=1}^{d_1}\left| \psi_p\right|^{-\frac{N+r_3}{2}}\psi_p^{-\alpha_{0 \psi_p}-1} \exp \Bigg\{-\frac{1}{2}\sum_{p=1}^{d_1}\Bigg(\sum_{i=1}^n \sum_{j=1}^J \psi^{-1}_p\Bigg(y_{i j p}-\boldsymbol{\beta}_p \boldsymbol{X}_{i j}-\boldsymbol{\Lambda}_p \boldsymbol{\omega}_{i j}\Bigg. \Bigg.\Bigg.
\Bigg.\Bigg.\Bigg.-\sum_{h=1}^{d_3} \boldsymbol{\phi}_{p h} \boldsymbol{Z}_{i j h}\Bigg)^2 \\
&+ \left(\boldsymbol{\Lambda}_p-\boldsymbol{\Lambda}_{0 p}\right) \left(\psi_p \boldsymbol{\Sigma}_{0 \Lambda_p}\right)^{-1}\left(\boldsymbol{\Lambda}_p-\boldsymbol{\Lambda}_{0 \Lambda_p}\right)^{\top}\Bigg)-\sum_{p=1}^{d_1} \beta_{0 \psi_p}\psi^{-1}_p\Bigg\},
\end{aligned}
\end{equation}
and  $p(\boldsymbol{\phi} \mid  \boldsymbol{\Lambda}, \boldsymbol{\Psi},\boldsymbol{ \beta} ,\boldsymbol{Y}, \boldsymbol{X}, \boldsymbol{Z})$  is proportional to
\begin{equation}\label{eq: CD_Phi}
\begin{aligned} 
&\left(2 \pi \right)^{-\frac{(N+K_n)d_1}{2}}\prod_{p=1}^{d_1}\left(\left| \psi_p\right|^{-\frac{N}{2}} \prod_{h=1}^{d_3}\left|\theta_{p h}^2 \bm{U}^{-1}\right|^{-\frac{1}{2}}\right)  \\
& \times \exp \Bigg\{-\frac{1}{2} \sum_{p=1}^{d_1} \Bigg(\sum_{i=1}^n \sum_{j=1}^J \psi^{-1}_p\Bigg(y_{i j p}-\boldsymbol{\beta}_p \boldsymbol{X}_{i j}-\boldsymbol{\Lambda}_p \boldsymbol{\omega}_{i j}-\sum_{h=1}^{d_3} \bm{\phi}_{p h} \boldsymbol{Z}_{i j h}\Bigg)^2 +\sum_{h=1}^{d_3}\bm{\phi}_{p h}\left(\theta_{p h}^2 \bm{U}^{-1}\right)^{-1} \bm{\phi}_{p h}^{\top}\Bigg)\Bigg\}.
\end{aligned}
\end{equation}
The precision matrices have conjugate Wishart full conditionals: 
\begin{equation}\label{eq: CD_T_xi}
\begin{aligned}
p\left(\boldsymbol{T}_{\xi}^{-1} \mid \boldsymbol{F}_2, \boldsymbol{A}_1, \boldsymbol{\alpha}_1\right) \stackrel{D}{=} \text { Wishart }_{r_2}\left(N+\mu_{0 \xi},\boldsymbol{R}_\xi\right),
\end{aligned}
\end{equation}
\begin{equation}\label{eq: CD_T_delta}
\begin{aligned}
&p\left(\boldsymbol{T}_\delta^{-1} \mid \boldsymbol{F}_2, \boldsymbol{X}, \boldsymbol{\Gamma}, \boldsymbol{A}_2, \boldsymbol{\alpha}_2\right)\stackrel{D}{=} \text { Wishart }_{r_1}\left(N+\mu_{0 \delta}, \boldsymbol{R}_\delta\right) ,
\end{aligned}
\end{equation}
where $\boldsymbol{R}_\xi=\left(\boldsymbol{R}_{0 \xi}^{-1}+\sum_{i=1}^n \boldsymbol{F}_{i(2)} \boldsymbol{C}_i\left(\boldsymbol{\alpha}_{i 1}\right)^{-1} \boldsymbol{F}_{i(2)}^\top\right)^{-1}$,  $\boldsymbol{R}_\delta=\left(\boldsymbol{R}_{0 \delta}^{-1}+\sum_{i=1}^n \boldsymbol{\Xi}_i \boldsymbol{G}_i^* \boldsymbol{D}_i\left(\boldsymbol{\alpha}_{i 2}\right)^{-1} \boldsymbol{G}_i^{*^\top} \boldsymbol{\Xi}_i^\top\right)^{-1}$.

 It follows from $p(\boldsymbol{B} \mid  \boldsymbol{\Gamma}, \mathbf{\Upsilon}) \propto p(\boldsymbol{\Gamma} \mid \boldsymbol{B}, \mathbf{\Upsilon}) p(\boldsymbol{B})$ that the conditional distribution of fixed effects $\boldsymbol{B}$ given $( \boldsymbol{\Gamma}, \mathbf{\Upsilon})$ follow a multivariate normal:
\begin{equation}\label{eq: CD_B}
\begin{aligned}
p(\boldsymbol{B} \mid  \boldsymbol{\Gamma}, \boldsymbol{\Upsilon}) \stackrel{D}{=} N_{r_2}\left(\boldsymbol{\mu}_B^*, \boldsymbol{\Omega}_B^*\right) , 
\end{aligned}
\end{equation}
where  $\boldsymbol{\Omega}_B^*=(\boldsymbol{\Sigma}_{0 B}^{-1}+ \boldsymbol{\Upsilon} ^{-1})^{-1}$ and $\boldsymbol{\mu}_B^*=\boldsymbol{\Omega}_B^*(\Sigma_{0 B}^{-1} \boldsymbol{B}_0+ \boldsymbol{\Upsilon} ^{-1} \sum_{i=1}^{n}\boldsymbol{\Gamma}_i)$. 
For the random effects covariance $\boldsymbol{\Upsilon}$,
$$
\begin{aligned}
&p(\boldsymbol{\Upsilon} ^{-1}\mid  \boldsymbol{\Gamma}, \boldsymbol{B}) \stackrel{D}{=} \text { Wishart }_{r_2}\left(n+\mu_{0 \Upsilon},\left(\boldsymbol{{R}_{0 \Upsilon}}^{-1}+\sum_{i=1}^n\left(\boldsymbol{\Gamma}_i- \boldsymbol{B}\right)\left(\boldsymbol{\Gamma}_i- \boldsymbol{B}\right)^\top\right)^{-1}\right). 
\end{aligned}
$$
If $\mathbf{\Upsilon}=\operatorname{diag}(\Upsilon_1, \ldots, \Upsilon_{r_2})$, then we have
\begin{equation}\label{eq: CD_Upsilon}
\begin{aligned}
p\left(\Upsilon_l ^{-1}\mid  \boldsymbol{\Gamma}, \boldsymbol{B}\right) \stackrel{D}{=} \operatorname{Gamma}\left(n+\alpha_{0 \Upsilon}, \beta_{0 \Upsilon}+\frac{1}{2} \sum_{i=1}^n a_{i l}^2\right),
\end{aligned}
\end{equation}
where $a_{i l}$ is the $l$-th component of vector $(\boldsymbol{\Gamma}_i -\boldsymbol{B})^\top $ for $l=1, \ldots, r_2$.

By Eq. \eqref{eq: CAR}, the conditional distribution 
\begin{equation}\label{eq: CD_rho_xi}
\begin{aligned}
p\left(\rho_{i \xi} \mid \boldsymbol{F}_2\right)\propto\left(\prod_{i=1}^n\left|\boldsymbol{I}_{T_i}-\rho_{i \xi} \boldsymbol{H}_{i \xi}\right|\right)^{\frac{r_2}{2}} I\left\{\rho_{i \xi} \in\left(\rho_{i \xi}^L, \rho_{i \xi}^U\right)\right\}
\times\exp \left\{-\frac{1}{2} \operatorname{tr}\left(\sum_{i=1}^n \boldsymbol{F}_{i(2)}^\top \boldsymbol{T}_{\xi}^{-1} \boldsymbol{F}_{i(2)}\left(\boldsymbol{I}_{T_i}-\rho_{i \xi} \boldsymbol{H}_{i \xi}\right)\right)\right\},
\end{aligned}
\end{equation}
where $\boldsymbol{I}_{T_i}$ is a $T_i$-dimensional identify matrix. 
Similarly, 
\begin{equation}\label{eq: CD_rho_delta}
\begin{aligned}
p(\rho_{i \delta} \mid \boldsymbol{F})\propto\left(\prod_{i=1}^n\left|\boldsymbol{I}_{T_i}-\rho_{i \delta} \boldsymbol{H}_{i \delta}\right|\right)^{\frac{r_1}{2}} I\left\{\rho_{i \delta} \in\left(\rho_{i \delta}^L, \rho_{i \delta}^U\right)\right\}
\times\exp \left\{-\frac{1}{2} \operatorname{tr}\left(\boldsymbol{G}_i^{* T} \boldsymbol{\Xi}_i^\top \boldsymbol{T}_\delta^{-1} \boldsymbol{\Xi}_i \boldsymbol{G}_i^*\left(\boldsymbol{I}_{T_i}-\rho_{i \delta} \boldsymbol{H}_{i \delta}\right)\right)\right\}  .
\end{aligned}
\end{equation}
The conditional distributions for the remaining parameters are provided in the Appendix \ref{sec:CD_SI}.  Among these,  the conditional distributions in Eqs.  \eqref{eq: CD_Beta-p}, \eqref{eq: CD_Lambda-p}, \eqref{eq: CD_T_xi}, \eqref{eq: CD_T_delta}, \eqref{eq: CD_B},  \eqref{eq: CD_Upsilon}  are standard normal, Wishart, and Gamma distributions, simulating observations from these conditional distributions is straightforward and fast. However, the conditional distributions in Eqs. \eqref{eq: CD_F}, 
\eqref{eq: CD_Psi}, \eqref{eq: CD_Phi} and \eqref{eq: CD_rho_xi}, \eqref{eq: CD_rho_delta} are non-standard and complex distributions, making direct sampling impractical. For these cases, we employ the Metropolis-Hastings (MH) algorithm. Detailed implementation of the MH algorithm is provided in Appendix \ref{sec:CD-MH_SI}.

\section{Simulation Studies}\label{sec:SS} 

We conducted comprehensive simulations to evaluate parameter recovery, computational efficiency, and prior sensitivity. 

In this simulation study, $\left\{y_{i j p}\right\}(i=1, \ldots, 30; j=1, \ldots, 6; $ $ p=1, \ldots, 9)$ are generated from the PFDSEM defined in Eqs. \eqref{eq: PFSEM}, \eqref{eq: structural_equation}, \eqref{eq: structural_heterogeneity}, \eqref{eq: covariance_matrices} and \eqref{eq: CAR}  with nine manifest variables that are related to three latent factors $\boldsymbol{\eta}_{i j}=\eta_{i j}$ and $\boldsymbol{\xi}_{i j}=(\xi_{i j 1}, \xi_{i j 2})^\top$, $\boldsymbol{\omega}_{i j}=(\eta_{i j},\xi_{i j 1}, \xi_{i j 2})^\top$. $\boldsymbol{Y}_{i j}$ is generated from a normal distribution $N_p(\boldsymbol{\beta}\boldsymbol{X}_{i j}+\boldsymbol{\Lambda}\boldsymbol{\omega}_{i j} +\int_{\mathcal{T}} \boldsymbol{C}(t)\boldsymbol{Z}_{i j}(t)d t, \boldsymbol{\Psi})$, where $\boldsymbol{X}_{i j}=(1,\boldsymbol{X}_{i j}^*)^\top=(1,{x}_{i j 1},{x}_{i j 2},{x}_{i j 3})^\top$, and $\boldsymbol{X}_{i j}^*$ is generated from a normal distribution $\boldsymbol{N_3}(\boldsymbol{0},\boldsymbol{I}_3)$, $\boldsymbol{Z}_{i j}(t)=({z}_{i j 1}(t),{z}_{i j 2}(t),{z}_{i j 3}(t))^\top$, $\boldsymbol{\Psi}=diag(\psi_1,\cdots,\psi_9)$.
The specifications of $\boldsymbol{\beta}$, $\boldsymbol{\Lambda}$, $\boldsymbol{C(t)}$ in relation to the measurement equation are given by 
$$\boldsymbol{\beta}=\left(\begin{array}{cccc}
\beta_{1,0}	& \beta_{1,1}	& \beta_{1,2}	& \beta_{1,3}\\
\beta_{2,0}	& \beta_{2,1}	& \beta_{2,2}	& \beta_{2,3}\\
...	&...	&...	&...	\\
\beta_{9,0}	& \beta_{P,1}	& \beta_{9,2}	& \beta_{9,3}\\
\end{array}\right), \quad
\boldsymbol{C(t)}=\left(\begin{array}{cccc}
c_{1 1}(t)	& c_{1 2}(t)	& c_{1 3}(t)\\
c_{2 1}(t)	& c_{2 2}(t)	& c_{2 3}(t)\\
...	&...	&...	\\
c_{9 1}(t)	& c_{9 2}(t)	& c_{9 3}(t)\\
\end{array}\right),
$$
$$\boldsymbol{\Lambda}^\top=\left(\begin{array}{ccccccccc}
1.0^*	&  \lambda_{21}	 & \lambda_{31}	& 0.0^*	& 0.0^*	 & 0.0^*& 0.0^*	 &0.0^*	 & 0.0^* \\
 0.0^*	& 0.0^*	& 0.0^*	 &1.0^*	 & \lambda_{52}	 & \lambda_{62}	& 0.0^*	& 0.0^*	& 0.0^* \\
 0.0^*	& 0.0^*	& 0.0^*	& 0.0^*	& 0.0^*	& 0.0^*	 &1.0^*	 & \lambda_{83}	 & \lambda_{93} 
\end{array}\right),
$$ where the values $1.0^*$ and $0.0^*$ are fixed for achieving an identified model.  $z_{i j h}(t)$ is generated  from a Gaussian process with  $z_{i j h}(t)\stackrel{D}{=}N(\{\sin (2 \pi t)+1.25\} / 2,\sigma_z^2)$ for $t \in[0,1]$, for the coefficient function, we set $c_{p h}(t)\stackrel{D}{=}N(\{\sin (2 \pi t)+1.25\} / 2,\sigma_c^2)$, where $\sigma_z^2=1$ and $\sigma_c^2=0.25$. 

The structural equation is defined by $\eta_{i j}=\gamma_{i 1} \xi_{i j 1}+\gamma_{i 2} \xi_{i j 2}+\delta_{i j}$, where $\xi_{i j 1}$ and $\xi_{i j 2}$ are respectively generated from a normal distribution $N(0,0.8)$. Since $\boldsymbol{\xi}_{i j}$ is two-dimensional and $\delta_{i j}$ is one-dimensional, we consider the following structures for $\boldsymbol{A}_1$ and $\boldsymbol{A}_2$ in modeling the between-variable covariances of $\boldsymbol{\xi}_{i j}$ and $\delta_{i j}$, respectively:
$$
\boldsymbol{A}_1=\left(\begin{array}{ll}
a_{11} & a_{12} \\
0.0 & a_{22}
\end{array}\right), \quad A_2=a_2.
$$
Here, $a_2=1.0$ is fixed for identification purpose.  For $i=1, \ldots, n, \boldsymbol{H}_{i \xi}=\boldsymbol{H}_{i \delta}=(h_{j l})_{J \times J}$ in which $h_{j l}=0$ with the exception of $h_{j, j+1}=h_{j+1, j}=1(j=1, \ldots, J-1)$, indicating that only the correlations between one-adjacent-time points, say time $j$ and time $j+1$, are considered.  For simplicity, we assume that $\rho_{i \xi}=\rho_{\xi}$ and $\rho_{i \delta}=\rho_\delta$ are invariant across individuals. In modeling the random structural coefficients with Eq. \eqref{eq: structural_heterogeneity}, we take $\boldsymbol{O}_i=\boldsymbol{I}_2, \boldsymbol{B}=(B_1, B_2)^\top, U_i=1$, and $\boldsymbol{\Upsilon}=\operatorname{diag}(\Upsilon_1,  \Upsilon_2)$.

 The true values of the unknown parameters are given by $\boldsymbol{\beta_p}=$ $(\beta_{p 0}, \beta_{p 1},\beta_{p 2},\beta_{p 3})$ $=(0.5, 0.8$ $, 0.8, 0.8)$, $ p=1, \ldots, 9 $, $ \lambda_{21}=\lambda_{31}=\lambda_{52}=\lambda_{62}=\lambda_{83}=\lambda_{93}=0.8$, $B_1=B_2=0.8$, $\Upsilon_1=\Upsilon_2=0.5$, $\rho_{\xi}=\rho_\delta=0.25$, $a_{11}=a_{22}=1.0$, and $a_{12}=0.5$, $\psi_p=1(p=1,2,3)$  and $\psi_p=0.8(p=4, \cdots, 9)$.  

We consider the following two cases: (\romannumeral1): $n=30$ with $J_1=\cdots=J_{30}=6$; (\romannumeral2): $n=80$ with $J_1=\cdots=J_{80}=6$. To investigate the sensitivity of bayesian estimates to inputs of hyperparameters in prior distributions, we consider the following two types of hyperparameters. 

\noindent  \textit{Type \textbf{\uppercase\expandafter{\romannumeral1}}}: $
 \boldsymbol{\mu}_0, 
 \boldsymbol{\Lambda}_{0}, \boldsymbol{B}_0$ are taken to be their corresponding true values, $\alpha_{0 \psi p}=10.0$ and $\beta_{0 \psi p}=8.0$ for $p=1, \ldots, 9$, $\Sigma_{0 {\beta_p}}, \Sigma_{0 \Lambda p}, \Sigma_{0 B}$ are diagonal matrices with diagonal elements $0.25, 
\alpha_{0 \Upsilon}=10.0, \beta_{0 \Upsilon}=5.0,
 \mu_{0 \xi}=10.0$, $\alpha_{0 \theta p}=1, \beta_{0 \theta p}=0.1$,and 
    $$\boldsymbol{R}_{0 \xi}=\left(\begin{array}{ll}
5.0 & 0.5 \\
0.5 & 5.0
\end{array}\right)
$$

\noindent \textit{Type \textbf{\uppercase\expandafter{\romannumeral2}}}: $\boldsymbol{\beta}_0, \boldsymbol{\Lambda}_{0}, \boldsymbol{B}_0$ are taken as $0.5, \boldsymbol{\Sigma}_{0 {\beta_p}}, \Sigma_{0 \Lambda p}$ and $\boldsymbol{\Sigma}_{0 B}$, are diagonal matrices with diagonal elements 10.0,and $\boldsymbol{R}_{0 \xi}=\operatorname{diag}(1.0, 1.0)$,
other hyperparameters are taken to be the same as those given in Type I.

In the implementation of the MH algorithm, we took $ \sigma_{\xi}^2=1.9$, $\sigma_\psi^2=2.0$,  $\sigma_{\rho_{\xi}}^2=20.1, \sigma_{\rho _{\delta}}^2=25.3$, and $\sigma_{\phi}^2=10.3$ in their corresponding proposal distributions to achieve approximate acceptance rates $0.34,0.32,0.33,0.34$ and $0.33$ , respectively.

\begin{figure}[!ht]
	\centering
	\begin{subfigure}{0.95\linewidth}
		\centering
        \includegraphics[width=0.5\textwidth]{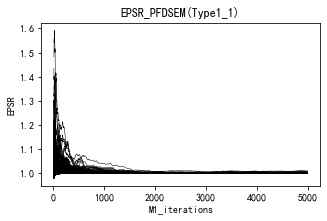}\includegraphics[width=0.5\textwidth]{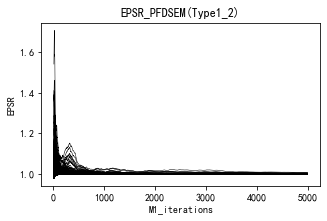}
	\end{subfigure}
	\centering
	\begin{subfigure}{0.95\linewidth}
		\centering
        \includegraphics[width=0.5\textwidth]{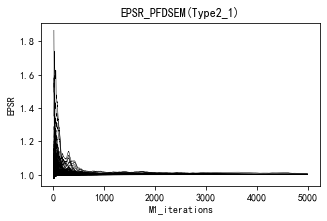}\includegraphics[width=0.5\textwidth]{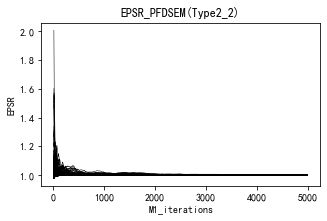}
	\end{subfigure}
    \caption{EPSR values of all parameters against iteration numbers in the simulation study}
    \label{fig:EPSR sim}
\end{figure}

Based on the above settings and the generated datasets, the proposed MCMC algorithm and variational inference  are used to get the Bayesian estimates on the basis of $100$ replications. In each of $100$ replications, $1000$ observations were collected to calculate the Bayesian estimates of unknown parameters after discarding $4000$ burn-in iterations.  To investigate the convergence of the algorithm, three parallel sequences of observations were generated from different starting values of the unknown parameters, and the EPSR values of all unknown parameters were computed and displayed in Fig. \ref{fig:EPSR sim}, 

The bias (Bias), standard deviation (SD) and the root mean square error (RMSE) between the Bayesian estimates and their true values are reported in Table \ref{tab: Performance-sim}, Fig. \ref{fig:est for phi}, and Fig. \ref{fig: est for phi1}, \ref{fig: est for phi2} and \ref{fig: est for phi3} in Appendix \ref{sec:SS_SI}. 

\setlength{\tabcolsep}{3pt} 
\begin{longtable}{cccccccccccccc}
\caption{Performance of the Bayesian estimates in the simulation study.}\label{tab: Performance-sim} \\ 
\hline \multirow{3}{*}{ par. } &  \multicolumn{6}{c}{ n=30, J=6 } &  & \multicolumn{6}{c}{ n=80, J=6 }  \\
\cdashline{2-7}[0.5pt/0.5pt]   \cdashline{9-14}[0.5pt/0.5pt] &  \multicolumn{3}{c}{ type1 } & \multicolumn{3}{c}{ type2 } & & \multicolumn{3}{c}{ type1 } & \multicolumn{3}{c}{ type2 } \\
 \cdashline{2-4}[0.5pt/0.5pt] \cdashline{5-7}[0.5pt/0.5pt]  \cdashline{9-11}[0.5pt/0.5pt] \cdashline{12-14}[0.5pt/0.5pt] 
& Bias & RMSE &SD & Bias & RMSE & SD &&  Bias & RMSE & SD & Bias & RMSE & SD \\
 \hline
\endfirsthead
\multicolumn{14}{c}%
{Table \thetable\ continued.} \\
\hline
\multirow{3}{*}{ par. } &  \multicolumn{6}{c}{ I=30,J=6 } &  & \multicolumn{6}{c}{ I=80,J=6 }  \\
\cdashline{2-7}[0.5pt/0.5pt]   \cdashline{9-14}[0.5pt/0.5pt]
  &  \multicolumn{3}{c}{ type1 } & \multicolumn{3}{c}{ type2 } & & \multicolumn{3}{c}{ type1 } & \multicolumn{3}{c}{ type2 } \\
 \cdashline{2-4}[0.5pt/0.5pt] \cdashline{5-7}[0.5pt/0.5pt]  \cdashline{9-11}[0.5pt/0.5pt] \cdashline{12-14}[0.5pt/0.5pt] 
& Bias & RMSE &SD & Bias & RMSE & SD &&  Bias & RMSE & SD & Bias & RMSE & SD \\
 \hline
\endhead
\hline \multicolumn{14}{c}{{Continued on next page}} \\ \hline
\endfoot
\hline
\endlastfoot
$\beta_{1 1}$&0.079 &0.187 &0.168 &0.038 &0.114 &0.107 & &0.096 &0.139 &0.101 &-0.004 &0.084 &0.084 \\
$\beta_{1 2}$&0.002 &0.141 &0.142 &-0.001 &0.101 &0.102 & &0.013 &0.105 &0.105 &0.012 &0.089 &0.089 \\
$\beta_{1 3}$&-0.003 &0.184 &0.182 &-0.009 &0.129 &0.126 & &-0.006 &0.098 &0.099 &-0.002 &0.083 &0.084 \\
$\beta_{1 4}$&0.024 &0.157 &0.157 &0.013 &0.106 &0.106 & &0.004 &0.097 &0.098 &0.004 &0.082 &0.082 \\
$\beta_{2 1}$&0.090 &0.171 &0.145 &0.054 &0.114 &0.100 & &0.089 &0.131 &0.098 &0.007 &0.083 &0.083 \\
$\beta_{2 2}$&0.003 &0.114 &0.115 &0.001 &0.084 &0.084 & &0.002 &0.084 &0.085 &0.001 &0.073 &0.073 \\
$\beta_{2 3}$&-0.005 &0.141 &0.142 &-0.010 &0.099 &0.099 & &-0.006 &0.088 &0.088 &-0.003 &0.077 &0.077 \\
$\beta_{2 4}$&0.013 &0.149 &0.150 &0.005 &0.111 &0.112 & &0.008 &0.086 &0.085 &0.008 &0.075 &0.075 \\
$\beta_{3 1}$&0.086 &0.166 &0.141 &0.050 &0.109 &0.096 & &0.084 &0.117 &0.082 &0.001 &0.069 &0.069 \\
$\beta_{3 2}$&-0.005 &0.119 &0.120 &-0.008 &0.089 &0.089 & &0.001 &0.078 &0.078 &0.001 &0.066 &0.066 \\
$\beta_{3 3}$&-0.017 &0.157 &0.157 &-0.020 &0.115 &0.114 & &0.001 &0.090 &0.090 &0.002 &0.079 &0.080 \\
$\beta_{3 4}$&0.010 &0.144 &0.143 &0.002 &0.105 &0.104 & &0.010 &0.089 &0.089 &0.009 &0.077 &0.077 \\
$\beta_{4 1}$&0.092 &0.140 &0.105 &0.067 &0.106 &0.082 & &0.094 &0.115 &0.065 &0.039 &0.070 &0.055 \\
$\beta_{4 2}$&0.009 &0.092 &0.092 &0.006 &0.078 &0.078 & &0.001 &0.054 &0.055 &0.001 &0.049 &0.049 \\
$\beta_{4 3}$&0.001 &0.114 &0.115 &-0.003 &0.096 &0.095 & &0.003 &0.059 &0.059 &0.004 &0.055 &0.055 \\
$\beta_{4 4}$&0.008 &0.114 &0.115 &0.003 &0.094 &0.095 & &-0.004 &0.058 &0.059 &-0.004 &0.054 &0.054 \\
$\beta_{5 1}$&0.099 &0.141 &0.101 &0.077 &0.112 &0.082 & &0.076 &0.103 &0.070 &0.029 &0.071 &0.064 \\
$\beta_{5 2}$&0.006 &0.084 &0.085 &0.004 &0.075 &0.076 & &-0.002 &0.057 &0.057 &-0.002 &0.052 &0.052 \\
$\beta_{5 3}$&0.008 &0.091 &0.091 &0.005 &0.076 &0.077 & &-0.001 &0.057 &0.057 &0.001 &0.054 &0.054 \\
$\beta_{5 4}$&0.001 &0.108 &0.109 &-0.003 &0.094 &0.095 & &-0.004 &0.053 &0.054 &-0.004 &0.049 &0.050 \\
$\beta_{6 1}$&0.080 &0.124 &0.095 &0.060 &0.098 &0.078 & &0.085 &0.108 &0.066 &0.038 &0.072 &0.060 \\
$\beta_{6 2}$&0.008 &0.085 &0.085 &0.005 &0.075 &0.075 & &-0.002 &0.055 &0.056 &-0.002 &0.051 &0.051 \\
$\beta_{6 3}$&-0.001 &0.097 &0.097 &-0.004 &0.082 &0.082 & &-0.003 &0.051 &0.052 &-0.002 &0.048 &0.049 \\
$\beta_{6 4}$&0.001 &0.090 &0.091 &-0.005 &0.078 &0.078 & &-0.006 &0.064 &0.065 &-0.005 &0.060 &0.061 \\
$\beta_{7 1}$&0.084 &0.136 &0.108 &0.062 &0.108 &0.089 & &0.085 &0.107 &0.066 &0.033 &0.070 &0.061 \\
$\beta_{7 2}$&0.001 &0.110 &0.110 &-0.001 &0.095 &0.095 & &-0.005 &0.065 &0.065 &-0.006 &0.060 &0.061 \\
$\beta_{7 3}$&0.002 &0.105 &0.106 &-0.001 &0.090 &0.091 & &-0.010 &0.063 &0.063 &-0.008 &0.060 &0.060 \\
$\beta_{7 4}$&-0.003 &0.103 &0.104 &-0.006 &0.089 &0.089 & &0.001 &0.058 &0.059 &0.001 &0.054 &0.055 \\
$\beta_{8 1}$&0.083 &0.115 &0.079 &0.064 &0.093 &0.067 & &0.074 &0.099 &0.066 &0.030 &0.068 &0.061 \\
$\beta_{8 2}$&0.010 &0.090 &0.089 &0.008 &0.079 &0.077 & &0.007 &0.055 &0.055 &0.006 &0.052 &0.052 \\
$\beta_{8 3}$&0.006 &0.108 &0.106 &0.004 &0.094 &0.093 & &-0.008 &0.048 &0.048 &-0.006 &0.046 &0.046 \\
$\beta_{8 4}$&-0.001 &0.093 &0.094 &-0.005 &0.082 &0.083 & &0.001 &0.054 &0.054 &0.001 &0.051 &0.052 \\
$\beta_{9 1}$&0.091 &0.123 &0.082 &0.071 &0.099 &0.069 & &0.079 &0.100 &0.061 &0.037 &0.067 &0.057 \\
$\beta_{9 2}$&-0.006 &0.101 &0.100 &-0.007 &0.089 &0.088 & &0.001 &0.056 &0.056 &0.001 &0.052 &0.053 \\
$\beta_{9 3}$&0.011 &0.080 &0.079 &0.008 &0.068 &0.068 & &-0.006 &0.059 &0.059 &-0.005 &0.056 &0.056 \\
$\beta_{9 4}$&0.007 &0.093 &0.092 &0.003 &0.083 &0.083 & &-0.002 &0.058 &0.059 &-0.002 &0.056 &0.056 \\
$\lambda_{1}$&-0.002 &0.050 &0.050 &0.005 &0.048 &0.048 & &-0.002 &0.031 &0.031 &0.001 &0.030 &0.030 \\
$\lambda_{2}$&-0.012 &0.053 &0.052 &-0.005 &0.050 &0.050 & &-0.005 &0.030 &0.030 &-0.003 &0.030 &0.030 \\
$\lambda_{3}$&-0.068 &0.119 &0.097 &-0.012 &0.075 &0.073 & &-0.021 &0.054 &0.049 &-0.003 &0.046 &0.046 \\
$\lambda_{4}$&-0.066 &0.123 &0.105 &-0.009 &0.080 &0.080 & &-0.020 &0.060 &0.057 &-0.001 &0.051 &0.052 \\
$\lambda_{5}$&-0.055 &0.124 &0.112 &0.008 &0.087 &0.087 & &-0.040 &0.078 &0.066 &-0.014 &0.060 &0.058 \\
$\lambda_{6}$&-0.056 &0.121 &0.108 &0.006 &0.083 &0.084 & &-0.029 &0.071 &0.065 &-0.004 &0.057 &0.057 \\
$\psi_{1}$&0.021 &0.136 &0.136 &0.030 &0.135 &0.132 & &-0.006 &0.106 &0.107 &-0.002 &0.106 &0.106 \\
$\psi_{2}$&0.007 &0.134 &0.134 &-0.003 &0.132 &0.132 & &0.005 &0.085 &0.085 &0.001 &0.083 &0.084 \\
$\psi_{3}$&0.026 &0.134 &0.131 &0.014 &0.129 &0.128 & &0.012 &0.082 &0.081 &0.008 &0.082 &0.082 \\
$\psi_{4}$&-0.074 &0.142 &0.122 &-0.047 &0.124 &0.116 & &-0.027 &0.076 &0.072 &-0.015 &0.071 &0.070 \\
$\psi_{5}$&-0.007 &0.097 &0.097 &-0.027 &0.099 &0.096 & &-0.007 &0.066 &0.066 &-0.015 &0.066 &0.065 \\
$\psi_{6}$&-0.013 &0.098 &0.098 &-0.035 &0.102 &0.097 & &0.001 &0.062 &0.062 &-0.009 &0.061 &0.060 \\
$\psi_{7}$&-0.039 &0.108 &0.101 &-0.012 &0.097 &0.097 & &-0.034 &0.081 &0.074 &-0.020 &0.074 &0.072 \\
$\psi_{8}$&0.001 &0.085 &0.085 &-0.022 &0.087 &0.084 & &0.005 &0.076 &0.076 &-0.004 &0.075 &0.075 \\
$\psi_{9}$&-0.007 &0.093 &0.094 &-0.027 &0.097 &0.094 & &0.001 &0.070 &0.070 &-0.007 &0.068 &0.069 \\
$\Upsilon_{1}$&-0.064 &0.119 &0.101 &-0.057 &0.113 &0.099 & &-0.041 &0.090 &0.080 &-0.036 &0.085 &0.078 \\
$\Upsilon_{2}$&-0.073 &0.114 &0.088 &-0.070 &0.113 &0.089 & &-0.057 &0.111 &0.095 &-0.049 &0.104 &0.093 \\
$\rho_\xi$&-0.006 &0.029 &0.029 &0.024 &0.029 &0.017 & &0.012 &0.035 &0.034 &0.031 &0.039 &0.024 \\
$\rho_\delta$&-0.038 &0.062 &0.049 &-0.032 &0.057 &0.047 & &-0.001 &0.043 &0.043 &0.001 &0.042 &0.042 \\
$\mathrm{T}_{\xi_{1 1}}$&0.009 &0.248 &0.250 &-0.010 &0.226 &0.228 & &-0.077 &0.150 &0.130 &-0.081 &0.147 &0.124 \\
$\mathrm{T}_{\xi_{1 2}}$&0.004 &0.124 &0.125 &-0.036 &0.124 &0.120 & &-0.025 &0.081 &0.077 &-0.040 &0.084 &0.074 \\
$\mathrm{T}_{\xi_{2 1}}$&0.004 &0.124 &0.125 &-0.036 &0.124 &0.120 & &-0.025 &0.081 &0.077 &-0.040 &0.084 &0.074 \\
$\mathrm{T}_{\xi_{2 2}}$&-0.005 &0.191 &0.193 &-0.018 &0.176 &0.176 & &-0.039 &0.130 &0.124 &-0.047 &0.125 &0.116 \\
$\mathrm{B}_{1}$&-0.042 &0.243 &0.242 &0.054 &0.160 &0.152 & &-0.014 &0.146 &0.147 &0.021 &0.124 &0.123 \\
$\mathrm{B}_{2}$&-0.095 &0.275 &0.260 &0.033 &0.152 &0.150 & &-0.040 &0.157 &0.154 &0.009 &0.121 &0.121 \\
\hline
$	\mathrm{Sum.}	$&	0.219	&	7.512	&	7.078	&	0.244	&	6.104	&	5.82	&	&	0.27	&	4.911	&	4.479	&	-0.115	&	4.229	&	4.089	\\
$	\mathrm{Mean}	$&	0.004	&	0.123	&	0.116	&	0.004	&	0.100	&	0.095	&	&	0.004	&	0.081	&	0.073	&	-0.002	&	0.069	&	0.067	\\ 																											\hline Time(min.)  &  \multicolumn{3}{c}{ 5.032 } & \multicolumn{3}{c}{ 4.425 } & & \multicolumn{3}{c}{ 12.807 } & \multicolumn{3}{c}{ 13.579 } \\
\end{longtable}
\setlength{\tabcolsep}{6pt} 

\begin{figure}[!ht]
	\centering
		\includegraphics[width=1\textwidth]{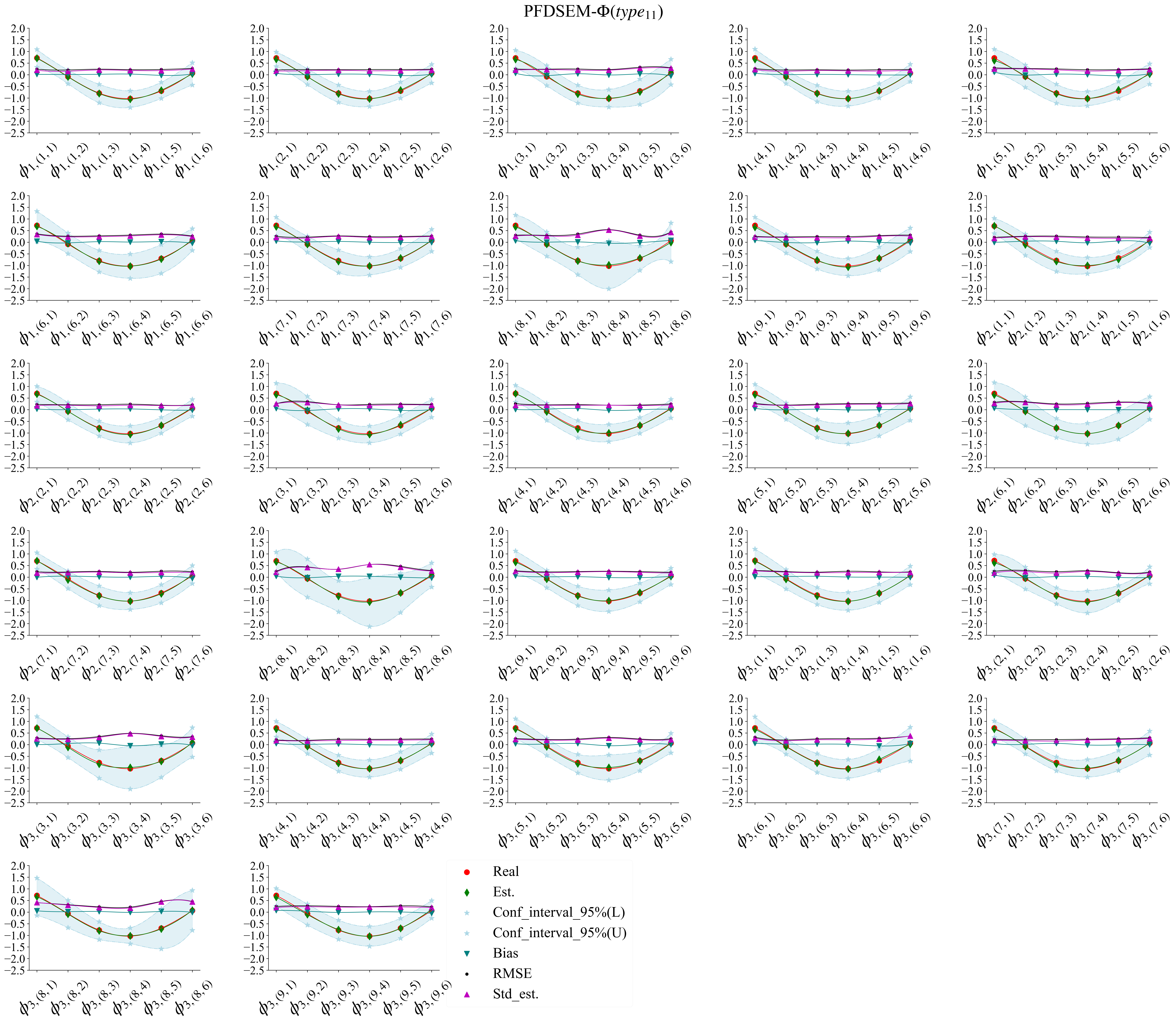}
	\caption{Bayesian estimation of functional parameters $\Phi$ ($type_{1 1}$) }
    \label{fig:est for phi}
\end{figure}

These results demonstrate that the proposed Bayesian estimation procedure performs well under the conditions examined. Across both sample sizes and prior specifications, the estimates exhibit consistently small bias and root mean square error, indicating accurate recovery of the true parameter values. The similarity of results between the informative (Type I) and weakly informative (Type II) prior settings suggests that the estimates are robust to the choice of prior hyperparameters, with no evidence of undue prior influence. The method also performs reliably with sample sizes as low as $n=30$, and estimation accuracy improves as expected when the sample size increases to $n=80$. The functional coefficient estimates in Fig. \ref{fig:est for phi} further confirm that the proposed approach successfully captures the underlying functional relationships with minimal bias and tight credible intervals.

\section{Applications of  Air Pollutant Emissions Data}\label{sec:Application} 

The objective of this section is to compile and develop a database of factors influencing air pollutant emissions, covering ten categories, and then to explore the extent of their impact on emission metrics. Our analysis framework is shown in Fig. \ref{fig:PFDSEM}. In our study, Spatial autocorrelation analysis (Moran's I), 
 spatial clustering analysis (Local Indicators of Spatial Association, LISA), 
and the proposed 
PFDSEMs will be utilized to investigate the spatio-temporal dynamic structural relationships between air pollutant indicators and their influencing factors.

\subsection{Study Area and Data}\label{sec:Data} 

Our study covered China's provinces, municipalities, and autonomous regions from January $1$, $2015$, to December $31$, $2020$. 
And our study incorporated the China provincial CO$_2$ emission inventory from the  China Emission Accounts and Datasets (CEADs) \citep{
xu2024china} 
and the emissions of nine air pollutants in china from the Multi-scale Emission Inventory of China (MEIC)\citep{
geng2024efficacy} 
as response variables for air pollutant emissions. 

Building upon prior research \citep{ghosh2010examining, jayanthakumaran2012co2, ozcan2013nexus, ZHU2021118323, LUO2022119242} and incorporating domain-specific characteristics of regional emissions, we systematically collected foundational determinants across ten conceptual dimensions. These dimensions include economic development level, standard of living, population structure, urbanization level, forest fires and natural disasters, vegetation greening, road traffic, technological innovation, household size, education level, environmental awareness, and meteorological factors. The complete dataset comprises 118 indicator variables. Meteorological data were obtained from the ERA5 database, environmental concern data were sourced from Baidu, and all other data were collected from the China Statistical Yearbook. Collinearity diagnosis and random forest methods were employed to screen the indicator features. Following this screening procedure, 49 indicators were retained as factor indicators for modeling and analysis. Detailed descriptions of the specific indicators are provided in Table \ref{tab: Indicators} (Appendix \ref{sec:Data_SI}). We categorize the conditioning factors into ten groups and construct the PFDSEM to explore the spatio-temporal structural relationships between each conditioning factor and the atmospheric pollutant indicators.

\subsection{Spatial Autocorrelation Analysis}\label{sec:saa}

\begin{figure}[ht]%
\centering
\includegraphics[width=1\textwidth]{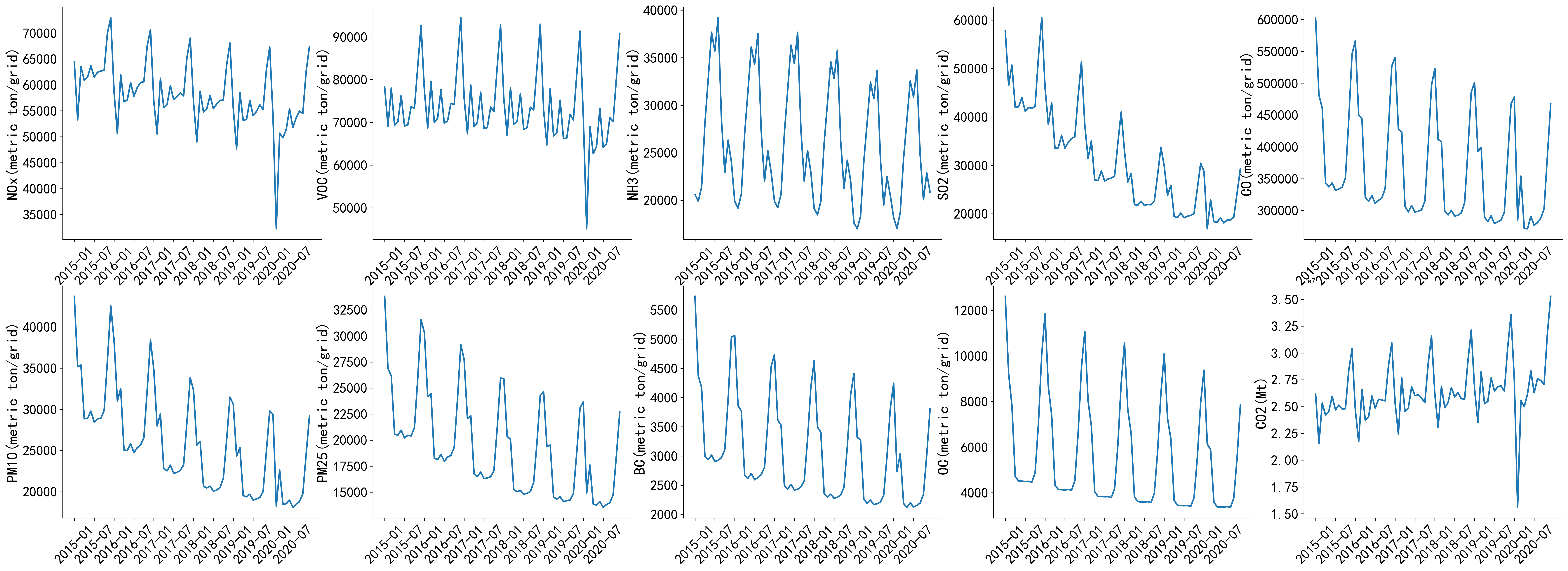}
\caption{Pollutant emission indicator time series}\label{PEITS}
\end{figure}

\begin{figure}[ht]%
\centering
\includegraphics[width=1\textwidth]{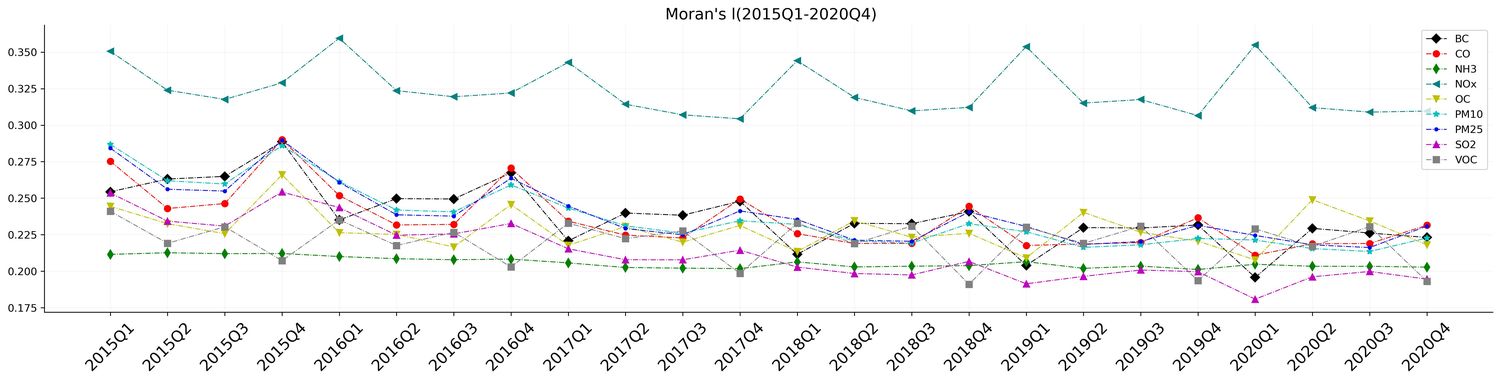}
\caption{ Moran's I values for each pollutant emissions in a geospatial weights matrix(quarter)}\label{fig:Moran's I}
\end{figure}

As depicted in Fig. \ref{PEITS}, spanning the interval from January $2015$ to December $2020$, a pronounced seasonal variation is observable in all evaluated air pollutant emission indices. The concentrations of SO$_2$, CO, PM$_{1 0}$, PM$_{2.5}$, and BC demonstrate a significant downward trajectory throughout the specified time frame, thereby affirming the effectiveness of china's strategies for carbon reduction in mitigating emissions. However, the time series merely delineate the shifting trends of various air pollutant emission indicator (APEI) over time, but failing to encapsulate the spatial characteristics inherent to the APEI. Accordingly, we have employed the global Moran's I index to articulate the spatial features of APEI across China and to evaluate its spatial correlation. Tables \ref{tab: Moran's I-PEI}(Appendix \ref{sec:SCA_SI}) and Fig. \ref{fig:Moran's I} present global Moran's I index 
 throughout the entirety of the study period. The significance levels at $1\%$ were consistently met across all seasons, signifying that the APEI exhibits a pronounced spatial clustering and positive spatial correlation in its distribution. A closer examination reveals that the global Moran indices for NO$_X$ and VOCs display no noteworthy declining or ascending trends, whereas the global Moran indices for the remaining indicators show a significant downward trend, indicating a gradual dilution of spatial agglomeration within APEI.

\begin{figure}[ht]%
 \centering
\includegraphics[width=1\textwidth]{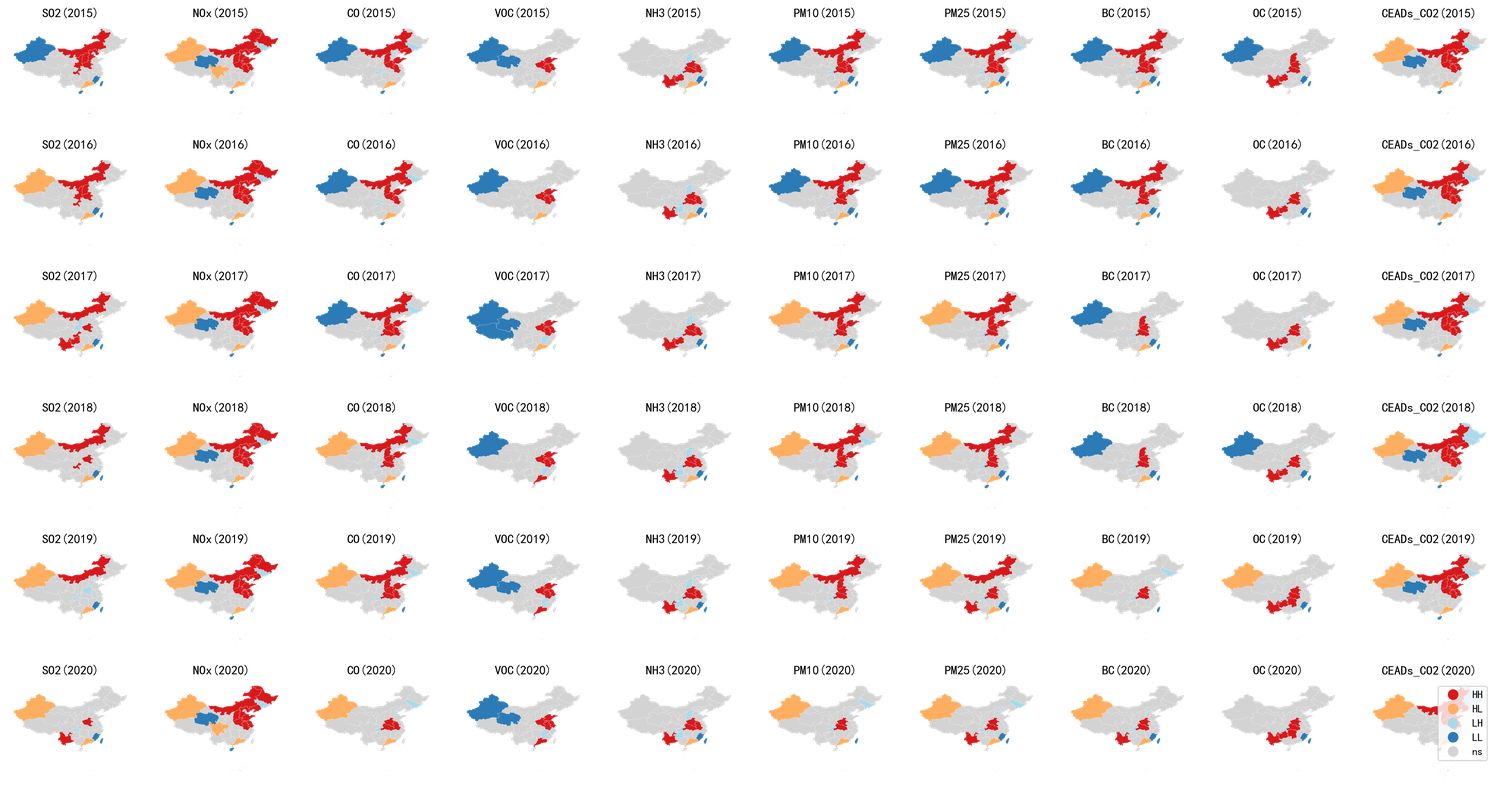}
\caption{ Lisa cluster of pollutant emission(2015-2020 year)}\label{fig:lisa year}
\end{figure}

Furthermore, adopting a spatio-temporal  analysis perspective of provincial disparities, this study engages in a comparative investigation of indicators spanning the years $2015$ to $2020$. The manuscript integrates the Local Indicators of Spatial Association (LISA)
 significance level test with Moran's I scatter plots to produce ``Moran's I Significance Maps'' adhering to a significance levels at $0.05$, as depicted in Fig. \ref{fig:lisa year}. And the analysis of the spatial autocorrelation patterns at the quarterly level is provided in Fig. \ref{fig41}- \ref{fig43} (Appendix \ref{sec:SCA_SI}). 

Upon thorough analysis, it has been determined that the divergent trends in pollutant emission indicators across various regions can be primarily attributed to several key factors: environmental policies and control measures, energy and industrial structure adjustments, regional economic development disparities, economic growth and industrial transformation, industrial and transportation development, agricultural activities and land use, and meteorological conditions changes. The detailed analysis results of the spatial agglomeration effects and spillover effects of various APEI indicators during the study period are  summarized in Table \ref{tab: Spatiotemporal correlation-APEI} (Appendix \ref{sec:SCA_SI}).

\subsection{Spatio-Temporal Dynamic Structural Relationships}\label{sec:STDSR}

Furthermore, we used the PFDSEM model to explore the spatio-temporal dynamic structural relationship between air pollutant emission indicators and their influencing factors. All observed variables were continuous variables, for simplicity, we assumed that they are continuous and 
$$
y_{i j p} \stackrel{D}{=}\left\{\begin{array}{l}
N\left(\boldsymbol{\beta_{p .}}\boldsymbol{X}_{i j}+\boldsymbol{\Lambda_{p .}}\boldsymbol{\omega}_{i j} +\sum_{k=1}^{K_n} \boldsymbol{\phi}_{p k} \boldsymbol{Z}_{i j k}, \psi_p\right), \quad p=1, \ldots, 10\\
N\left(\tilde{\beta_p}+\boldsymbol{\Lambda_{p .}}\boldsymbol{\omega}_{i j}, \psi_p\right), \quad p=11, \ldots, 54,
\end{array}\right.
$$ 
where the first $10$ variables $(Y_1,\cdots,Y_{1 0})$ were interpreted as ``Emission inventory'', while the variables $(Y_{1 1},\cdots, Y_{1 5})$, $(Y_{1 6},\cdots, Y_{2 1})$, $(Y_{2 2 },\cdots, Y_{2 6})$, $(Y_{2 7},\cdots, Y_{3 1})$, $(Y_{3 2},\cdots, Y_{3 7})$, $(Y_{3 8},\cdots, Y_{4 1})$, $(Y_{4 2},Y_{4 3}, Y_{4 4})$, $(Y_{4 5},Y_{4 6}, Y_{4 7})$, ($Y_{4 8}$, $Y_{4 9}$, $Y_{5 0}$),  ($Y_{5 1}, \cdots, Y_{5 4}$) were interpreted as ``Meteorology'', ``Economic'', ``Population structure'', ``Forest fires \& natural disasters'', ``Vegetation Cover\& Sci-Tech innovation'', ``Road \& Traffic'', ``Urbanization level'', ``Family scale'', ``Education'', and  ``Environmental Concern'' respectively. 
Scalar covariates $(X_1, X_2, X_3)$ corresponded to ``Environmental governance'', ``Government intervention'', and ``Foreign direct investment (FDI)'', respectively. Functional covariates $(Z_1(t), Z_2(t))$ represented ``Sea level pressure'' and ``2 m temperature''. To ensure scale uniformity, we standardized all raw data using complete observations.
Which corresponded to latent variables $\eta_{i j}$, $\xi_{i j 1}, \cdots,\xi_{i j 1 0}$. We choose the non-overlapping structure of latent factors according to the meaning of the questions. For the structural equation model (Eq. \eqref{eq: structural_equation}), we considered the following model:
$$
\begin{aligned}
 \eta_{i j}=\sum_{r_2=1}^{10}\gamma_{i r_2} \xi_{i j r_2}+\delta_{i j} .
\end{aligned}
$$
 For the individual-level mixed effect in Eq. \eqref{eq: structural_heterogeneity}, we considered  $\boldsymbol{B}=(B_{1}, \cdots, B_{10})$, and $\boldsymbol{\Upsilon}=\operatorname{diag}(\Upsilon_1, \cdots,  \Upsilon_{10})$.  Because $\boldsymbol{\xi}_{i j}=(\xi_{i j 1}, \cdots, \xi_{i j 10})^\top$ was ten-dimensional and $\delta_{i j}$ was one-dimensional, we considered the following structures for $\boldsymbol{A}_1$ and $A_2$ in modeling the between-variable covariances of $\boldsymbol{\xi}_{i j}$ and $\delta_{i j}$, respectively: 
$$
\boldsymbol{A}_1=\left(\begin{array}{cccccccccc}
a_{1 1}	&	a_{1 2}	&	a_{1 3}	&	a_{1 4}	&	a_{1 5}	&	a_{1 6}	&	a_{1 7}	&	a_{1 8}	&	a_{1 9}	&	a_{1, 10}	\\
0	&	a_{2 2}	&	a_{2 3}	&	a_{2 4}	&	a_{2 5}	&	a_{2 6}	&	a_{2 7}	&	a_{2 8}	&	a_{2 9}	&	a_{2, 10}	\\
0	&	0	&	a_{3 3}	&	a_{3 4}	&	a_{3 5}	&	a_{3 6}	&	a_{3 7}	&	a_{3 8}	&	a_{3 9}	&	a_{3, 10}	\\
0	&	0	&	0	&	a_{4 4}	&	a_{4 5}	&	a_{4 6}	&	a_{4 7}	&	a_{4 8}	&	a_{4 9}	&	a_{4, 10}	\\
0	&	0	&	0	&	0	&	a_{5 5}	&	a_{5 6}	&	a_{5 7}	&	a_{5 8}	&	a_{5 9}	&	a_{5, 10}	\\
0	&	0	&	0	&	0	&	0	&	a_{6 6}	&	a_{6 7}	&	a_{6 8}	&	a_{6 9}	&	a_{6, 10}	\\
0	&	0	&	0	&	0	&	0	&	0	&	a_{7 7}	&	a_{7 8}	&	a_{7 9}	&	a_{7, 10}	\\
0	&	0	&	0	&	0	&	0	&	0	&	0	&	a_{8 8}	&	a_{8 9}	&	a_{8, 10}	\\
0	&	0	&	0	&	0	&	0	&	0	&	0	&	0	&	a_{9 9}	&	a_{9, 10}	\\
0	&	0	&	0	&	0	&	0	&	0	&	0	&	0	&	0	&	a_{10, 10}	\\	
\end{array}\right), \quad A_2=a_2 ,
$$
where $a_2=1$ was fixed for identification purpose. 

For $i=1, \ldots, 30, \boldsymbol{H}_{i \xi}=\boldsymbol{H}_{i \delta}=(h_{j l})_{6 \times 6}$ was a matrix containing the neighboring two positions' information whose elements were given by $h_{j l}=0$ with exception of $h_{j, j+1}=h_{j+1, j}=1$ for $j=1,2, \ldots, 5$. That is, we considered the correlated structures for $\boldsymbol{\xi}_{i j}$ and $\delta_{i j}$ involving one-adjacent-time points.
The specified prior distributions: 
$\alpha_{0 \psi p}=10.0$ and $\beta_{0 \psi p}=8.0$ for $p=1, \ldots, 54$; $\mu_{0 \Upsilon}=10.0, R_{0 \Upsilon}=5.0$; $\Sigma_{0 \beta}, \Sigma_{0 \Lambda r}, \Sigma_{0 B}$ were diagonal matrices with diagonal elements $0.25$, and  $\boldsymbol{R}_{0 \xi}$ were diagonal matrices with diagonal elements $1$; $\mu_{0 \xi}=20.0$, $\alpha_{0 \theta p}=1, \beta_{0 \theta p}=0.1$,  $\boldsymbol{\beta}_0, \boldsymbol{\Lambda}_{0}, \boldsymbol{B}_0$ were taken to be their corresponding bayesian estimates obtained via non-informative prior distributions.

\begin{figure}[ht]
\centering
\includegraphics[width=1\textwidth]{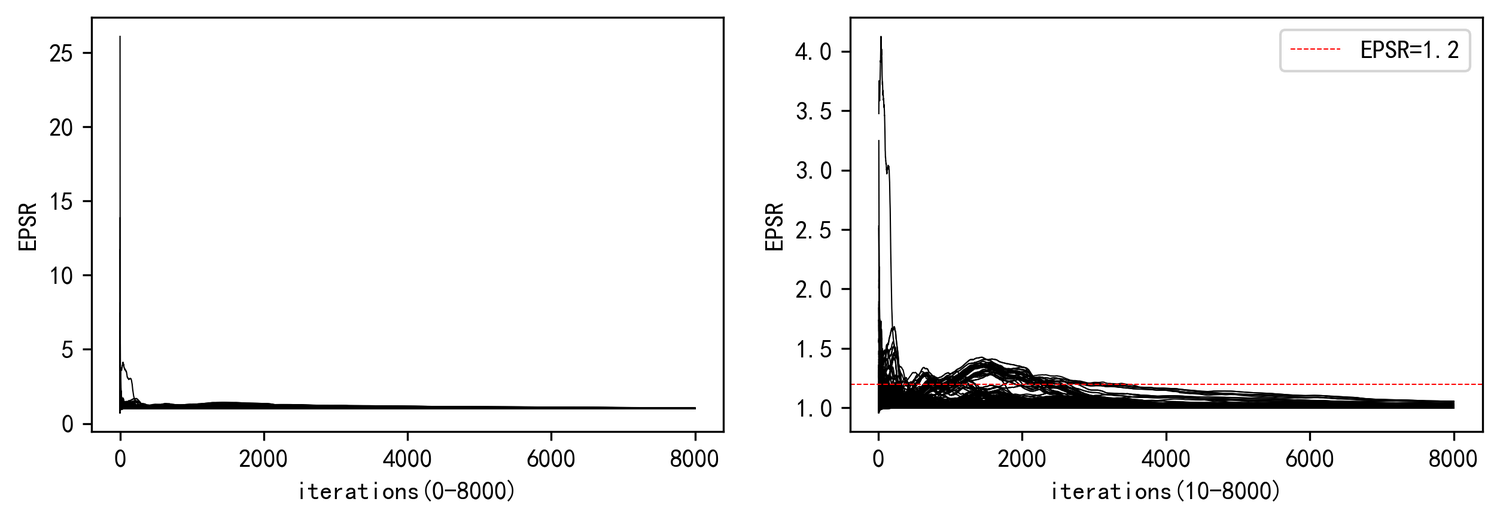}
\caption{EPSR of PFDSEM (real example) }\label{fig:EPSR real exam.}
	\end{figure}

\begin{figure}[ht]%
\centering
\includegraphics[width=1\textwidth]{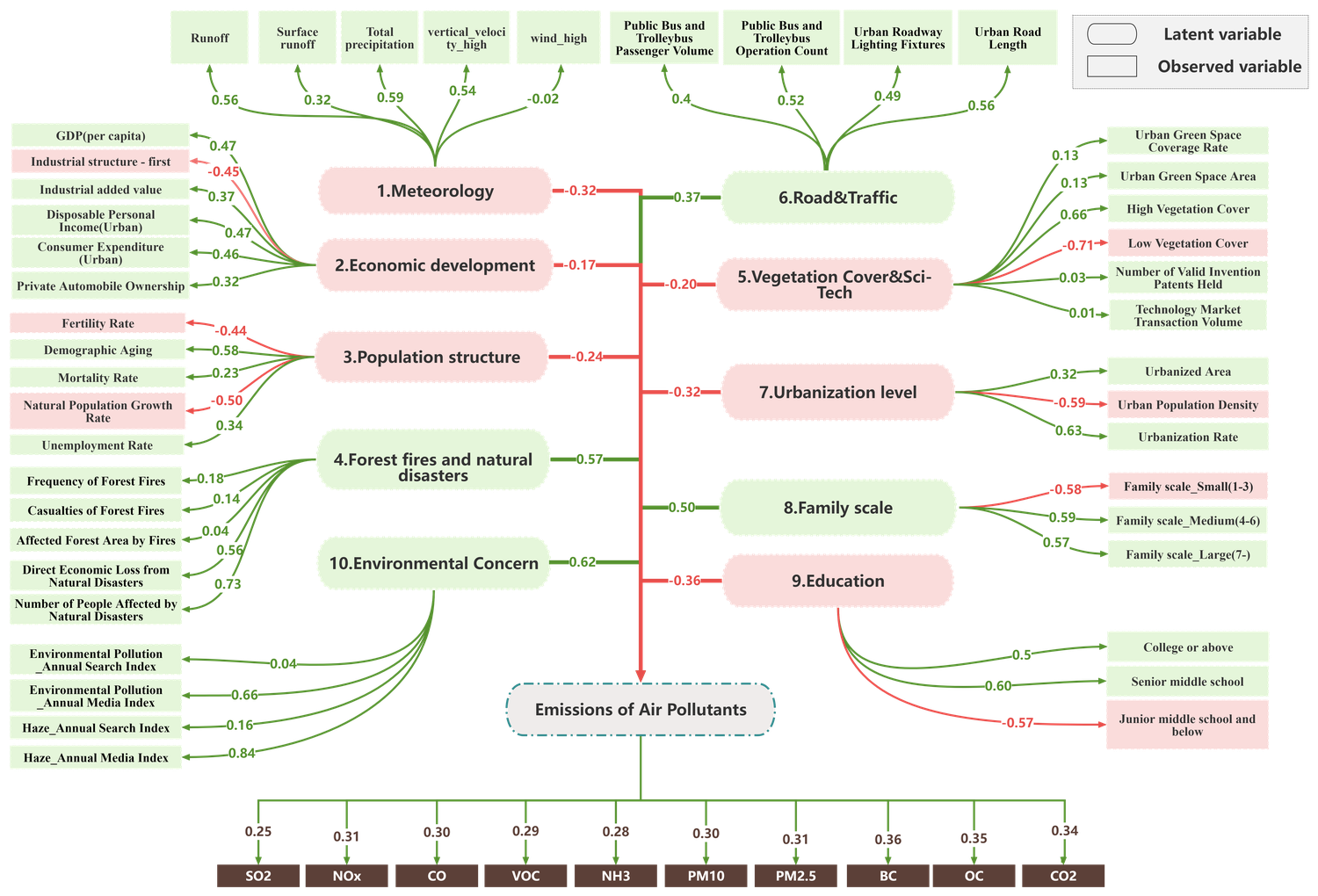}
\caption{Pollutant emissions  PFDSEM model}\label{fig:PFDSEM-APE}
\end{figure}

\begin{figure}[ht]%
\centering
\includegraphics[width=1\textwidth]{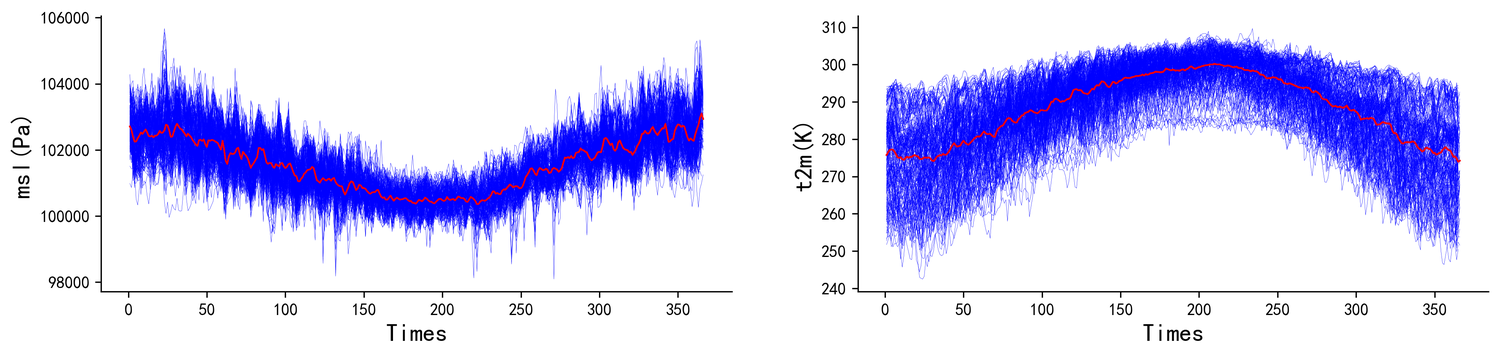}
\caption{Functional data time series}\label{fig:FDTS}
\end{figure}

\begin{figure}[ht]%
\centering
\includegraphics[width=1\textwidth]{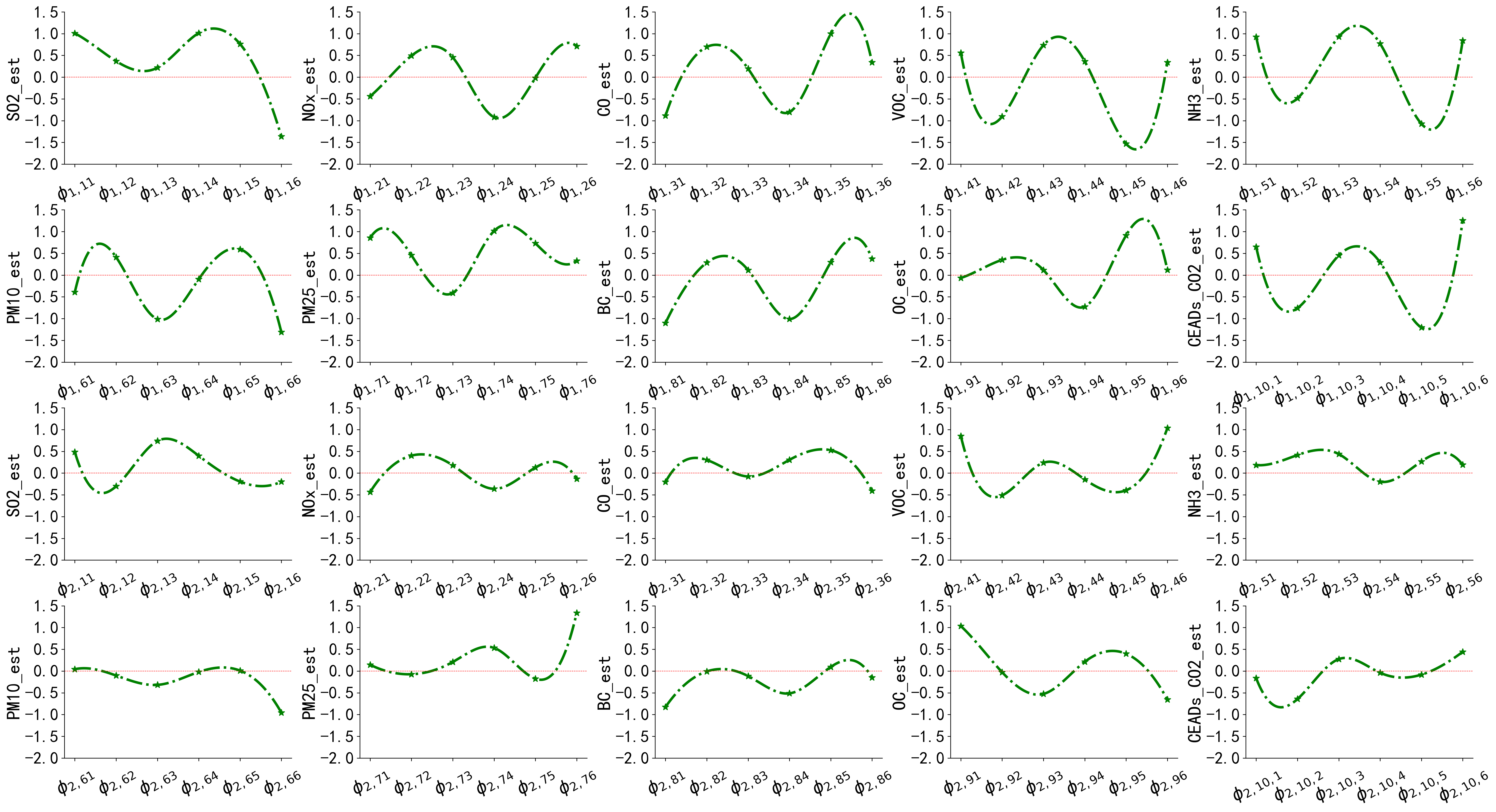}
\caption{Pollutant emissions  PFDSEM model (Estimated of functional coefficients,where phi1 and phi2 corresponds to the functional variable $Z_1(t)$ and $Z_2(t)$ respectively, that is, Sea level pressure and 2 m temperature)}\label{fig:est for phi (real)}
\end{figure}

\begin{figure}[!ht]%
\centering
\includegraphics[width=1\textwidth]{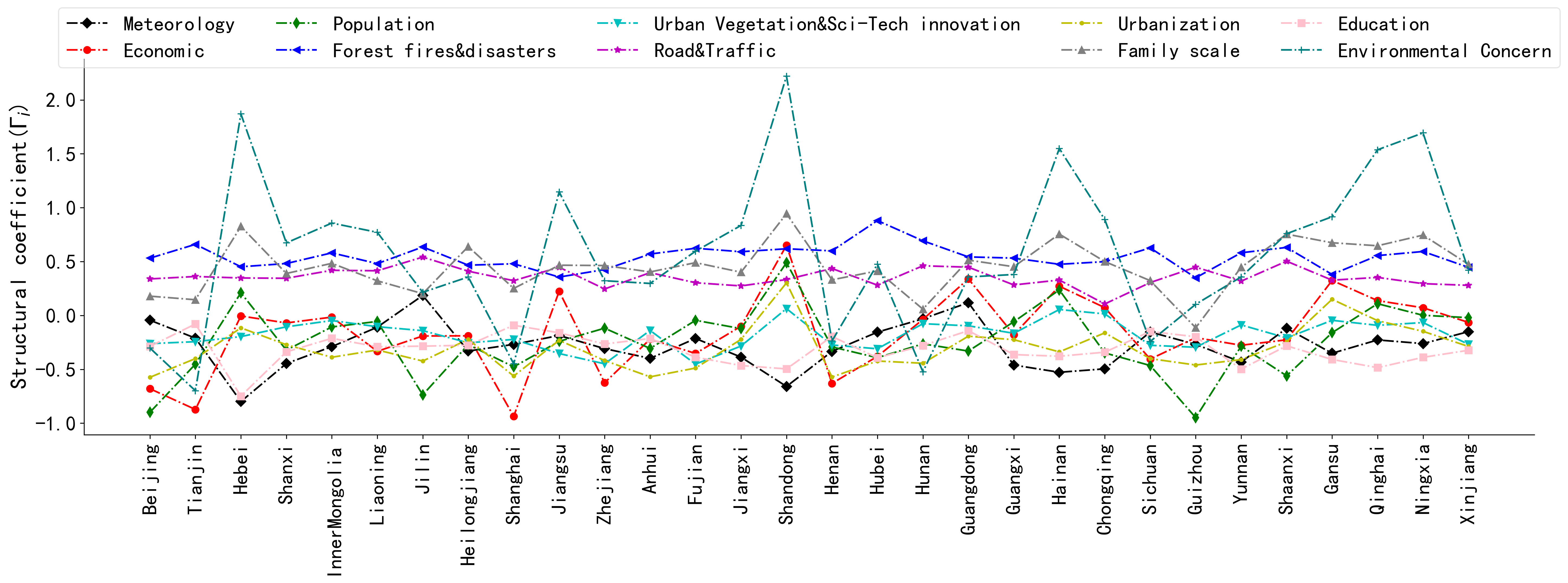}
\caption{Pollutant emissions  PFDSEM model(Structural coefficient)}\label{fig:Structural coefficient (real)}
\end{figure}

\begin{figure}[!ht]%
 \centering
\includegraphics[width=1\textwidth]{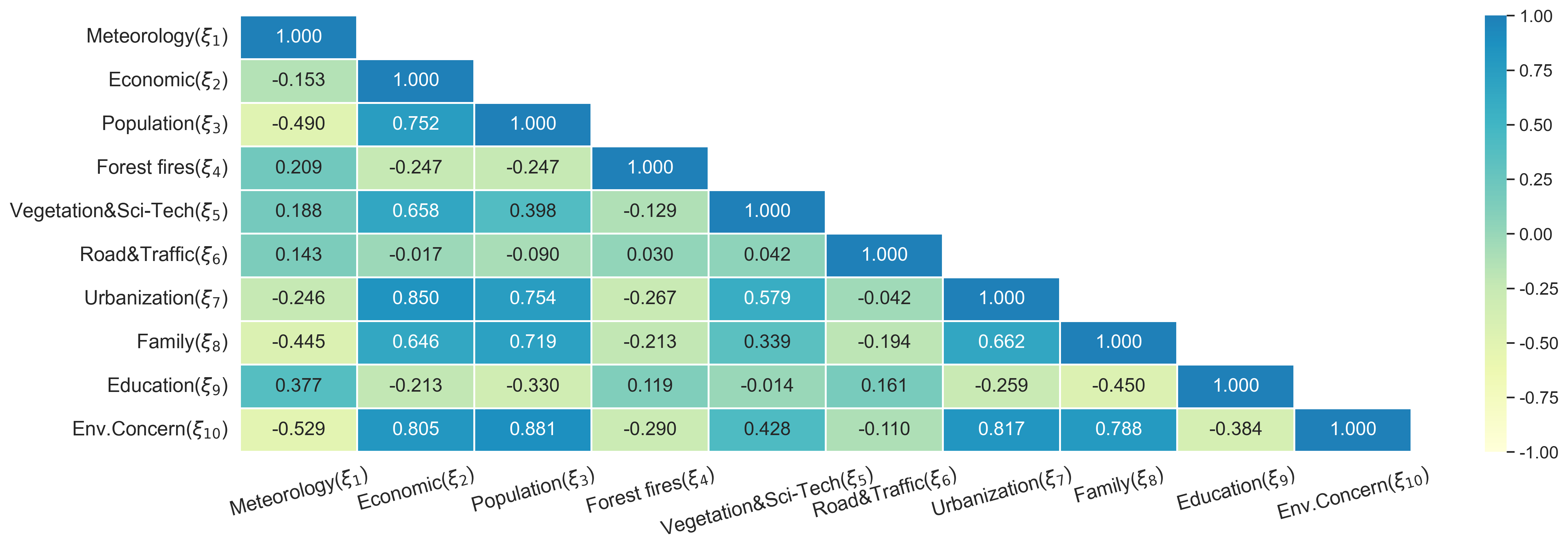}
\caption{Correlation coefficient of $\xi_1$ ... $\xi_{1 0}$}\label{fig:Correlation coefficient (real)}
\end{figure}

Here, we utilized the empirical Bayes method, which used not only the form of the likelihood, but also the observed data values to determine the prior because we had no prior knowledge about the unknown parameters. To investigate the convergence of the algorithm, three parallel sequences of observations were generated from different starting values of the unknown parameters, and the EPSR values of all unknown parameters were computed and displayed in Fig. \ref{fig:EPSR real exam.}. It showed that the estimated potential scale reduction (EPSR) are all below 1.2, signifying the convergence of the model. 
The Bayesian estimates and their standard error estimates of the unknown parameters in  model  were obtained with 8000 observations after 7000 burn-in iterations, the results were reported in Tables \ref{tab: fitting-APEI},  \ref{tab: Performance-APE} (Appendix \ref{sec:STDSR_SI}) and Fig. \ref{fig:PFDSEM-APE}-\ref{fig:Correlation coefficient (real)}, \ref{fig: fitting results (real)_1}-\ref{fig: fitting results (real)_5} (Appendix \ref{sec:STDSR_SI}) .

\begin{table}[htbp]
\centering
\caption{Model fitting effect of APEI.}\label{tab: fitting-APEI}
\begin{tabular}{ccccccccccc} 
\hline Ind. & SO2 & NOx  &CO& VOC & NH3&PM10& PM2.5 & BC&OC&CO2\\
\hline RMSE &0.092 & 0.100 & 0.062 & 0.151 & 0.170 & 0.064 & 0.067 & 0.095 & 0.147 & 0.133 \\
      SD& 0.131 & 0.157 & 0.156 & 0.149 & 0.151 & 0.156 & 0.163 & 0.186 & 0.182 & 0.178 \\
      $R^2$ & 0.685 & 0.776 & 0.896 & 0.583 & 0.421 & 0.891 & 0.887 & 0.839 & 0.675 & 0.714 \\
\hline
\end{tabular}
\end{table}
 
Based on the estimated individual-level parameters $\boldsymbol{B}$ and $\boldsymbol{\Upsilon}$ and their corresponding standard errors (Table \ref{tab: Performance-APE} and Fig. \ref{fig:PFDSEM-APE}), all influencing factors considered in this study exert statistically significant effects on the nine types of atmospheric pollutant emissions and CO$_2$ emissions. 
In particular,   Meteorological conditions$\left(\xi_1\right)$, Economic development$\left(\xi_2\right)$, Population structure$\left(\xi_3\right)$, Vegetation Cover and Sci-Tech innovation$\left(\xi_5\right)$,  Urbanization level$\left(\xi_7\right)$,  Education$\left(\xi_9\right)$ are identified as inhibitors of emissions indices. Conversely,   Forest fires and natural disasters$\left(\xi_4\right)$,  Road and Traffic$\left(\xi_6\right)$,  Family scale$\left(\xi_8\right)$,  Environmental Concern$\left(\xi_{1 0}\right)$ are found to exacerbate the intensity of emissions. Further comprehensive and meticulous analysis of the estimation results for the loadings ($\boldsymbol{\Lambda}$), scalar covariate coefficients ($ \boldsymbol{\beta}$), and functional covariate coefficients ($\boldsymbol{\Phi}$) described in Table \ref{tab: Performance-APE} and Fig. \ref{fig:PFDSEM-APE},\ref{fig:FDTS} and \ref{fig:est for phi (real)}. Detailed conclusions are provided in Appendix \ref{sec:SR_SI}.

From the estimated values $\hat{\boldsymbol{\Upsilon}}=(\hat{\Upsilon}_1,\hat{\Upsilon}_2,\ldots,\hat{\Upsilon}_{10})^\top=(0.50,0.53,0.60,0.51,0.47,0.50,0.51,0.56,0.53,0.84)$ and the structural coefficients $\boldsymbol{ \Gamma_i}=(\gamma_{i,1},\cdots,\gamma_{i,10})$  in Fig. \ref{fig:Structural coefficient (real)}, we observe that the structural relationships between atmospheric pollutant emission indicators and their driving factors exhibit notable heterogeneity across cities. 
A detailed discussion of the factors underlying these regional differences is provided in Appendix \ref{sec:STDD_SI}.

Furthermore, the Bayesian estimates of the adjacent-time correlation parameters $\hat{\rho}_{\xi}=0.545$ and $\hat{\rho}_{\delta}=0.555$, are significantly different from zero, as their posterior credible intervals exclude zero. This finding confirms that the structural relationships between air pollutant emission indicators and their driving factors exhibit substantial temporal dependence, and moreover, that these correlation patterns are non-stationary over the study period. This temporal dependence can be attributed to a combination of evolving socio-economic and environmental conditions, including changes in economic activity, shifts in energy consumption patterns, the implementation of policy interventions, and technological advancements. A detailed discussion of the factors underlying these temporal dependencies is provided in Appendix \ref{sec:STDD_SI}.

Building on the analysis of spatio-temporal heterogeneity, we now examine the interrelationships among the latent driving factors themselves. The estimated covariance matrix $\boldsymbol{T}_\xi$ reveals the correlation structure among the ten latent variables ($\xi_1,\ldots,\xi_{10}$), as visualized in Fig. \ref{fig:Correlation coefficient (real)}. Several key factors exhibit significant positive correlations, including economic development($\xi_2$), urbanization level($\xi_7$), population structure($\xi_3$), family scale($\xi_8$), environmental concern($\xi_{10}$), and technological innovation($\xi_5$). This pattern indicates that these factors do not operate independently; rather, they are mutually reinforcing within a complex socio-economic-environmental system. Such synergistic effects imply that policy interventions targeting one factor may indirectly influence others, underscoring the need for integrated policy approaches. A detailed analysis of these correlations and their substantive interpretations is provided in Appendix \ref{sec:Corr_SI}.

\section{Discussion}\label{sec:Disc}

In this study, we have proposed a Partially Functional Dynamic Structural Equation Model (PFDSEM) designed to address the methodological challenges posed by multi-resolution environmental data. By integrating functional covariates via basis expansion, scalar covariates within a standard regression framework, and latent variables within a dynamic SEM structure, the PFDSEM provides a coherent analytical platform for quantifying complex, temporally dynamic, and province-specific associations between atmospheric pollutant emissions and their multifaceted determinants.

\subsection{Key Environmental Findings and Their Implications}

Our analysis of 30 Chinese provinces over the period 2015--2020 reveals several findings with direct relevance to environmental policy and management:

\begin{enumerate}
    \item[i)]\textit{Factors associated with emissions.} Our analysis reveals a clear dichotomy in the direction of associations between the ten categories of socio-environmental factors and atmospheric pollutant emissions (Fig. \ref{fig:PFDSEM-APE}). On the one hand, meteorological conditions, economic development, vegetation cover, technological innovation, urbanization levels, and educational attainment all exhibit strong negative associations with emissions. These findings are consistent with the hypothesis that economic modernization, when accompanied by technological upgrading and improved environmental awareness, can contribute to emission reductions. On the other hand, forest fires and natural disasters, road and traffic conditions, household size, and environmental concern (measured as search intensity) exhibit positive associations with emission intensity. The positive association of environmental concern likely reflects reverse feedback: higher pollution levels trigger greater public attention, a pattern documented in the environmental health literature, while the positive association with road and traffic indicators underscores the persistent challenge of transportation-related emissions despite vehicle efficiency improvements. The strength of these associations (both negative and positive ) varies substantially across provinces, as evidenced by the estimated random effects variances $\hat{\Upsilon}_l$ and the structural coefficients $\boldsymbol{ \Gamma_i}=(\gamma_{i,1},\cdots,\gamma_{i,10})$ in Fig. \ref{fig:Structural coefficient (real)}, suggesting that the emission-reducing benefits of economic and technological development, as well as the emission-exacerbating pressures from infrastructure and natural hazards, are not uniformly realized and may depend on province-specific institutional and structural conditions.

 \item[iii)]\textit{Temporal dynamics and policy persistence.} The estimated temporal dependence parameters ($\hat{\rho}_{\xi}=0.545$, $\hat{\rho}_{\delta}=0.555$, both with 95\% credible intervals excluding zero) confirm that the structural relationships among latent factors exhibit substantial temporal persistence. 
 This finding has direct policy relevance: the effects of environmental regulations, economic incentives, or technological interventions on emissions are likely to be realized over extended time horizons, and short-term evaluations may substantially underestimate their cumulative impact.

 \item[iv)]\textit{Seasonal heterogeneity in meteorological effects.} The estimated functional coefficients for sea-level pressure and 2-meter temperature (Fig. \ref{fig:est for phi (real)}) reveal pronounced seasonal patterns in the strength of meteorological associations with emissions. These patterns reflect well-known atmospheric processes: winter temperature inversions trap pollutants near the surface in many Chinese cities, while summer atmospheric conditions favor dispersion. The ability to quantify these seasonal effects with full posterior uncertainty represents a significant advance over approaches that aggregate meteorological variables to seasonal averages, which can mask important sub-seasonal variability.

 \item[v)]\textit{Inter-factor correlations and policy complexity.} The estimated latent factor correlation matrix (Fig. \ref{fig:Correlation coefficient (real)}) reveals strong positive associations among several key factors, including economic development, urbanization, and technological innovation. These interdependencies imply that policy interventions targeting one factor may have cascading indirect associations with others, complicating the design of effective emission reduction strategies. For example, policies promoting urbanization without complementary investments in public transportation and green infrastructure may produce unintended increases in traffic-related emissions, offsetting potential efficiency gains from urban agglomeration.
\end{enumerate}

\subsection{Methodological Contributions and Comparison with Existing Approaches}

The PFDSEM extends the Generalized Random Coefficient Structural Equation Model \citep{song2012bayesian} in a non-trivial direction by accommodating functional covariates measured at higher temporal resolutions than the scalar outcomes and predictors. This extension is achieved through basis expansion, which projects the infinite-dimensional functional data onto a finite-dimensional space while preserving the continuity and smoothness of the underlying processes. Compared to ad-hoc approaches that discretize functional data or pre-aggregate to period averages, the basis expansion framework offers three advantages: (i) it retains within-period information that may be predictive of extreme events; (ii) it propagates uncertainty from the functional representation coherently through all model levels; and (iii) the Bayesian P-spline formulation provides a principled mechanism for controlling smoothness through data-driven penalty parameters.

The explicit separation of temporal dynamics (via the CAR structure) and inter-variable covariances (via the LMC Cholesky parameterization) enhances model interpretability. The CAR parameter $\rho$ directly quantifies the persistence of temporal dependence, while the LMC matrices $\boldsymbol{A}_1$ and $\boldsymbol{A}_2$ capture the contemporaneous correlation structure among latent factors. This separation is particularly valuable in environmental applications, where temporal persistence (e.g., policy inertia, technological lock-in) and structural correlations (e.g., economic-urbanization linkages) often reflect distinct underlying processes.

\subsection{Limitations and Future Research Directions}

Several limitations merit acknowledgment. First, despite controlling for complex structures, the observational design precludes strong causal claims; unmeasured confounders or reverse causality (e.g., pollution affecting migration or health) cannot be ruled out. Future research could integrate instrumental variables or natural experiments. Second, the CAR(1) structure, while parsimonious, may not capture all temporal dynamics; extensions to ARMA or stochastic volatility models could be explored. Third, the current treatment of functional covariates assumes additive effects; including function-on-function interactions could reveal more complex dependencies. Finally, computational demands remain substantial; variational Bayesian approximations or scalable MCMC algorithms (e.g., Hamiltonian Monte Carlo) warrant investigation. Despite these limitations, PFDSEM provides a principled, flexible framework for integrating multi-resolution data in environmental systems, and our empirical findings offer an evidence base for designing adaptive, seasonally-targeted pollution control policies. 

Despite these limitations, the PFDSEM provides a principled and flexible framework for integrating multi-resolution environmental data within a dynamic structural equation modeling paradigm. The empirical findings offer a quantitative evidence base for designing temporally adaptive and regionally differentiated environmental policies, and the methodological framework opens several avenues for further development at the intersection of functional data analysis, dynamic latent variable modeling, and environmental statistics.

\section*{Declarations}

\subsection*{Funding}
The authors were partially supported by National Natural Science Foundation
of China Grants (Nos. 12292980,12292983); And the topic `Robust Resampling Techniques for Big Data' of the National Key Research and Development Program(No.2022YFA1003701); And the Scientific Research Fund Project of Yunnan Education Department(No.2025Y0052), and Scientific Research and Innovation Project of Postgraduate Students in the Academic Degree of YunNan University (No.KC-24249954).
 
\subsection*{Conflict of interest/Competing interests}
The authors declare no competing interests. 

\bibliographystyle{apalike}
\bibliography{refs}

\begin{thebibliography}{}

\bibitem[Apergis and Payne, 2010]{apergis2010renewable}
Apergis, N. and Payne, J.~E. (2010).
\newblock Renewable energy consumption and economic growth: Evidence from a
  panel of oecd countries.
\newblock {\em Energy Policy}, 38(1):656--660.

\bibitem[Besag, 2018]{Besag1974SpatialIA}
Besag, J. (2018).
\newblock Spatial interaction and the statistical analysis of lattice systems.
\newblock {\em Journal of the Royal Statistical Society: Series B
  (Methodological)}, 36(2):192--225.

\bibitem[Chaudhary and Bisai, 2018]{chaudhary2018factors}
Chaudhary, R. and Bisai, S. (2018).
\newblock Factors influencing green purchase behavior of millennials in india.
\newblock {\em Management of Environmental Quality: An International Journal},
  29(5):798--812.

\bibitem[Gelman et~al., 2013]{gelman2013bayesian}
Gelman, A., Carlin, J.~B., Stern, H.~S., and et~al. (2013).
\newblock {\em Bayesian Data Analysis}.
\newblock Chapman and Hall/CRC, Boca Raton, FL, 3 edition.

\bibitem[Gelman and Hill, 2006]{gelman2006data}
Gelman, A. and Hill, J. (2006).
\newblock {\em Data Analysis Using Regression and Multilevel/Hierarchical
  Models}.
\newblock Analytical Methods for Social Research. Cambridge University Press,
  New York.

\bibitem[Geman and Geman, 1984]{Geman1984Gibbs}
Geman, S. and Geman, D. (1984).
\newblock Stochastic relaxation, gibbs distributions, and the bayesian
  restoration of images.
\newblock {\em IEEE Transactions on Pattern Analysis and Machine Intelligence},
  PAMI-6(6):721--741.

\bibitem[Geng et~al., 2024]{geng2024efficacy}
Geng, G., Liu, Y., Liu, Y., and et~al. (2024).
\newblock Efficacy of {China}'s clean air actions to tackle {PM}$_{2.5}$
  pollution between 2013 and 2020.
\newblock {\em Nature Geoscience}, 17(10):987--994.

\bibitem[Ghosh, 2010]{ghosh2010examining}
Ghosh, S. (2010).
\newblock Examining carbon emissions economic growth nexus for india: A
  multivariate cointegration approach.
\newblock {\em Energy Policy}, 38(6):3008--3014.

\bibitem[Hastie et~al., 2009]{Hastie2009ESL}
Hastie, T., Tibshirani, R., and Friedman, J. (2009).
\newblock {\em The Elements of Statistical Learning: Data Mining, Inference,
  and Prediction}.
\newblock Springer, New York, 2nd edition.

\bibitem[Jayanthakumaran et~al., 2012]{jayanthakumaran2012co2}
Jayanthakumaran, K., Verma, R., and Liu, Y. (2012).
\newblock Co2 emissions, energy consumption, trade and income: A comparative
  analysis of china and india.
\newblock {\em Energy Policy}, 42:450--460.

\bibitem[J\"{o}reskog, 1969]{Joreskog_1969}
J\"{o}reskog, K.~G. (1969).
\newblock A general approach to confirmatory maximum likelihood factor
  analysis.
\newblock {\em Psychometrika}, 34(2):183--202.

\bibitem[Klein and Kneib, 2016]{Klein2016ScaleDependentPriors}
Klein, N. and Kneib, T. (2016).
\newblock {Scale-Dependent Priors for Variance Parameters in Structured
  Additive Distributional Regression}.
\newblock {\em Bayesian Analysis}, 11(4):1071 -- 1106.

\bibitem[Lang and Brezger, 2004]{Lang2004BayesianPSplines}
Lang, S. and Brezger, A. (2004).
\newblock Bayesian p-splines.
\newblock {\em Journal of Computational and Graphical Statistics}, 13.

\bibitem[Lee, 2007]{Lee2007SEM}
Lee, S.-Y. (2007).
\newblock {\em Structural Equation Modeling: A Bayesian Approach}.
\newblock Wiley Series in Probability and Statistics. John Wiley \& Sons,
  Chichester, UK.

\bibitem[Luo et~al., 2022]{LUO2022119242}
Luo, S., Zhu, Y., and Chen, S.~X. (2022).
\newblock Episode based air quality assessment.
\newblock {\em Atmospheric Environment}, 285:119242.

\bibitem[Ozcan, 2013]{ozcan2013nexus}
Ozcan, B. (2013).
\newblock The nexus between carbon emissions, energy consumption and economic
  growth in middle east countries: A panel data analysis.
\newblock {\em Energy Policy}, 62:1138--1147.

\bibitem[Pourahmadi, 1999]{Pourahmadi1999}
Pourahmadi, M. (1999).
\newblock Joint mean-covariance models with applications to longitudinal data:
  unconstrained parameterisation.
\newblock {\em Biometrika}, 86(3):677--690.

\bibitem[Ramsay and Silverman, 2005]{Ramsay2005FDA}
Ramsay, J.~O. and Silverman, B.~W. (2005).
\newblock {\em Functional Data Analysis}.
\newblock Springer, New York, 2nd edition.

\bibitem[Rauf et~al., 2018]{rauf2018structural}
Rauf, A., Zhang, J., Li, J., and et~al. (2018).
\newblock Structural changes, energy consumption and carbon emissions in china:
  Empirical evidence from ardl bound testing model.
\newblock {\em Structural Change and Economic Dynamics}, 47:194--206.

\bibitem[Robert and Casella, 2004]{robert2004monte}
Robert, C.~P. and Casella, G. (2004).
\newblock {\em Monte Carlo Statistical Methods}.
\newblock Springer-Verlag, 2nd edition.

\bibitem[Song et~al., 2018]{song2018environmental}
Song, M., Peng, J., Wang, J., and et~al. (2018).
\newblock Environmental efficiency and economic growth of china: A ray
  slack-based model analysis.
\newblock {\em European Journal of Operational Research}, 269(1):51--63.

\bibitem[Song et~al., 2012]{song2012bayesian}
Song, X., Tang, N., and Chow, S. (2012).
\newblock A bayesian approach for generalized random coefficient structural
  equation models for longitudinal data with adjacent time effects.
\newblock {\em Computational Statistics \& Data Analysis}, 56(12):4190--4203.

\bibitem[Stephen et~al., 1998]{Brooks01121998}
Stephen, P., Brooks, Andrew, and Gelman (1998).
\newblock General methods for monitoring convergence of iterative simulations.
\newblock {\em Journal of Computational and Graphical Statistics}, 7:434--455.

\bibitem[Tanner et~al., 1987]{Tanner01061987}
Tanner, M.~A., Wong, W.~H., and et~al. (1987).
\newblock The calculation of posterior distributions by data augmentation.
\newblock {\em Journal of the American Statistical Association},
  82(398):528--540.

\bibitem[Wood, 2017]{Wood2017GAM}
Wood, S.~N. (2017).
\newblock {\em Generalized Additive Models: An Introduction with R}.
\newblock Chapman and Hall/CRC, New York, 2nd edition.

\bibitem[Xi, 2023]{xi2023china}
Xi, J. (2023).
\newblock Holding high the great banner of socialism with chinese
  characteristics to strive in unity for the great rejuvenation of the chinese
  nation and to build a modern socialist country in all respects -- report at
  the 20th national congress of the communist party of china.
\newblock {\em Dangshi Bocai}, 22:46.

\bibitem[Xu et~al., 2024]{xu2024china}
Xu, J., Guan, Y., Oldfield, J., and et~al. (2024).
\newblock China carbon emission accounts 2020-2021.
\newblock {\em Applied Energy}, 360:122837.

\bibitem[Zhu et~al., 2021]{ZHU2021118323}
Zhu, Y., Liang, Y., and Chen, S.~X. (2021).
\newblock Assessing local emission for air pollution via data experiments.
\newblock {\em Atmospheric Environment}, 252:118323.

\end{thebibliography}

\section*{Author Biography}

\begin{biography}{\includegraphics[width=0.14\textwidth]{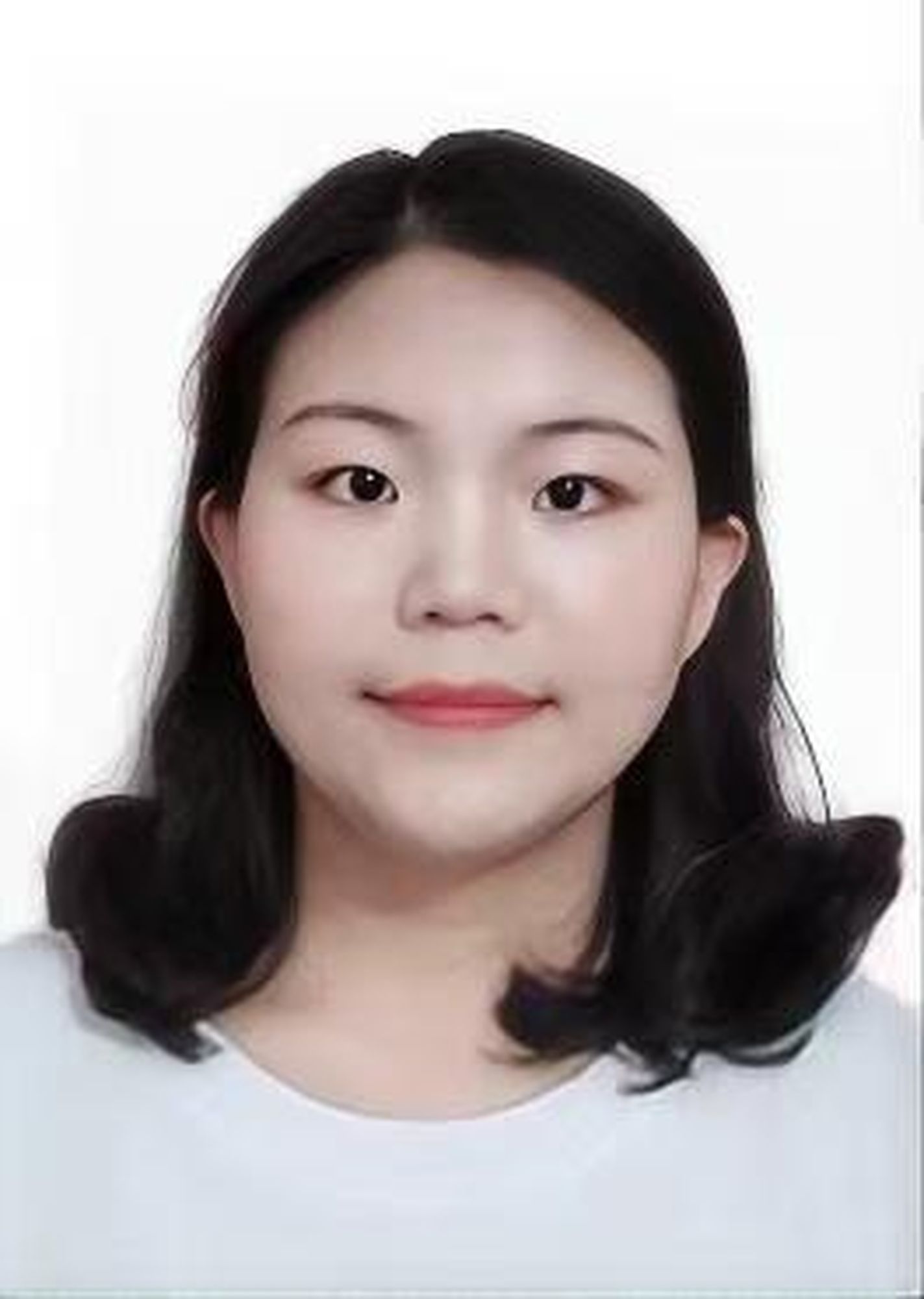}}{\textbf{Xinwen Liu\orcidlink{0000-0003-0563-9228}} is a Ph.D. student at the School of Mathematics and Statistics, Yunnan University (Southwest Joint Graduate School), supervised by Professor Tang Niansheng and Academician Song Xi Chen. She has been pursuing her Ph.D. at Yunnan University from 2023 to the present, and since March 2025, she has been a visiting scholar at the Center for Statistical Science, Peking University. She led one provincial‑level research project during the period 2024–2025. Her research interests include complex spatio‑temporal big data analysis, machine learning, Bayesian statistics, causal inference, deep learning, etc. }
\end{biography}

\begin{biography}{\includegraphics[width=0.14\textwidth]{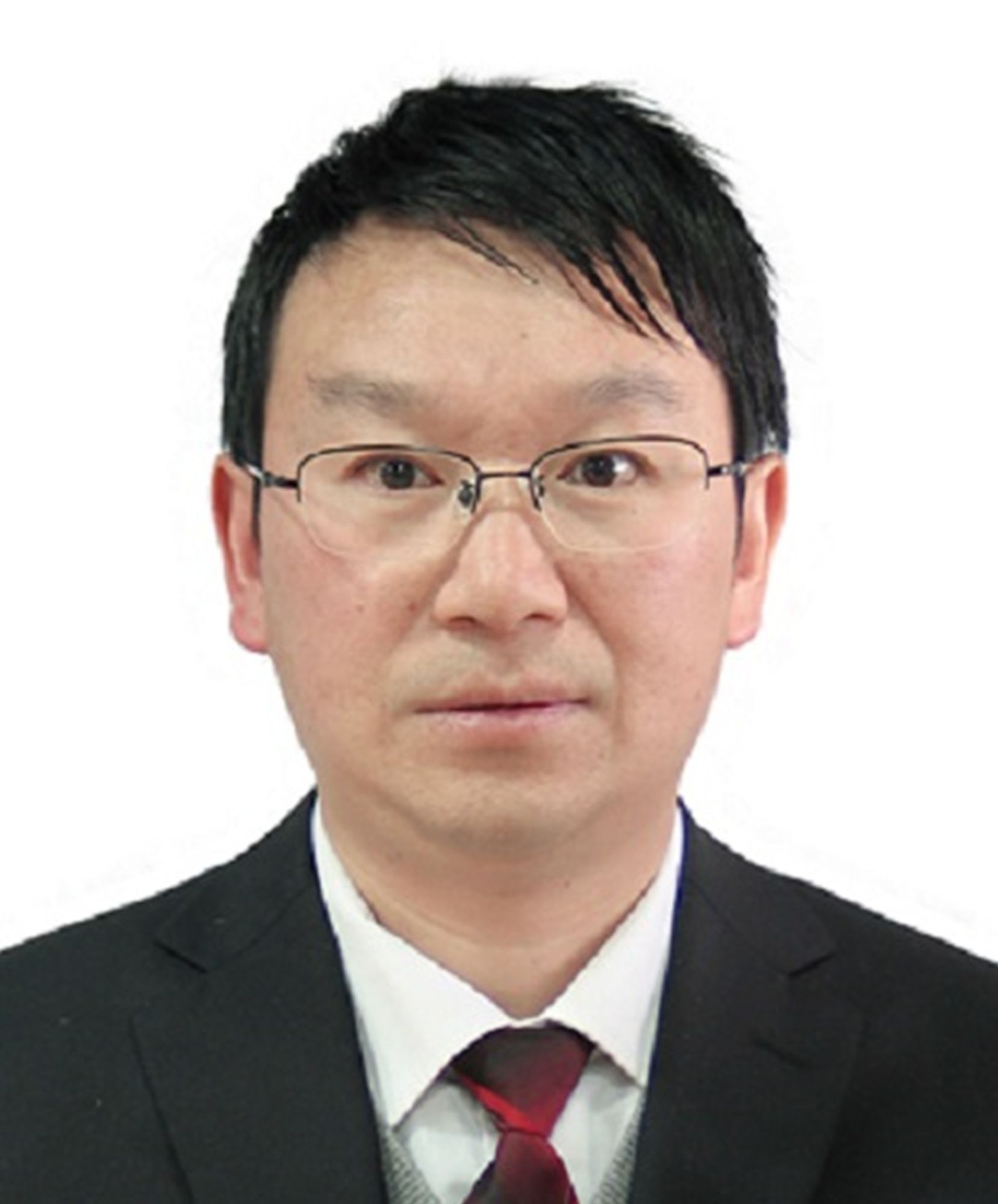}}{\textbf{Niansheng Tang \orcidlink{0000-0001-7033-3845}} is a professor at the School of Mathematics and Statistics of Yunnan University and doctoral supervisor, currently serves as Vice President of Yunnan University. He is a recipient of the ``National Outstanding Youth Science Foundation,'' a Ministry of Education Changjiang Scholar Distinguished Professor, Deputy Chairman of the Chinese Society of Applied Statistics, member of the Steering Committee for Undergraduate Statistics Programs in Higher Education Institutions of the Ministry of Education, awardee of the ``New Century Excellent Talents Support Plan'' by the Ministry of Education, a distinguished expert under the National Program for Special Support of Young and Middle-aged Outstanding Talents with Significant Contributions, enjoys a special allowance from the State Council, an elected member of the International Statistical Institute (ISI), a fellow of the Institute of Mathematical Statistics (IMS), chief scientist of national key research projects, and serves as associate editor and editorial board member for several internationally renowned statistical journals such as CIMS. His research areas include biostatistics, Bayesian statistics, statistical diagnostics, missing data analysis, high-dimensional data analysis, survival data analysis, among others. He has published over 200 academic papers in prestigious journals like Annals of Statistics, JRSSB (Journal of the Royal Statistical Society, Series B), JOE (Journal of Econometrics), JASA (Journal of the American Statistical Association), and Biometrika.}
\end{biography}

\begin{biography}{\includegraphics[width=0.14\textwidth]{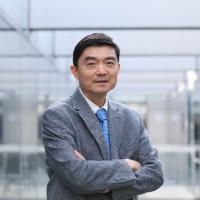}}{\textbf{Song Xi Chen \orcidlink{0000-0002-2338-0873}} is a Chair Professor in Statistics and Data Science at Tsinghua University. He is a Fellow of the Institute of Mathematical Statistics, the American Statistical Association, and the American Association for the Advancement of Science. He serves as President of the Society for Probability and Statistics of China and is a member of the Chinese Academy of Sciences. His research interests include ultra‑high dimensional statistical inference, statistical and machine learning, and statistical applications in earth science. He has published over 140 refereed papers.}
\end{biography}

\end{document}